\documentstyle[11pt]{article}

\newcommand{\A}{ {\cal A} }
\newcommand{\B}{ {\cal B} }
\newcommand{\C}{ \mbox{\bf C} }

\newcommand{\iz}{ i_{o} }
\newcommand{\jz}{ j_{o} }
\newcommand{\N}{ \mbox{\bf N} }
\newcommand{\R}{ \mbox{\bf R} }
\newcommand{\T}{ \mbox{\bf T} }
\newcommand{\X}{ {\cal X} }

\newcommand{\CC}{ {\cal C} }
\newcommand{\Z}{ \mbox{\bf Z} }

\newcommand{\dd}{ \underline{D} }
\newcommand{\hh}{ \underline{H} }
\newcommand{\sps}{ \underline{S} }
\newcommand{\xx}{ \underline{X} }
\newcommand{\hbblock}{ ( \hh_{B} \ ; \ B \mbox{ block of } \sigma ) }
\newcommand{\dhb}{ \dd \setminus \cup_{B}  \hh_{B} }
\newcommand{\cq}{ C_{Q} }
\newcommand{\crr}{ C_{R} }

\newcommand{\ee}{ \varepsilon }
\newcommand{\eealt}{ NC_{\ee - alt} (m) }
\newcommand{\hacc}{ NC_{h - acc} }

\newcommand{\aiqaik}{ a_{i_{1}} \cdots a_{i_{k}} }
\newcommand{\xiqxik}{ X_{i_{1}} \cdots X_{i_{k}} }
\newcommand{\ziqzik}{ z_{i_{1}} \cdots z_{i_{k}} }

\newcommand{\funuk}{ f_{1} , \ldots ,f_{k} }
\newcommand{\gunuk}{ g_{1} , \ldots ,g_{k} }
\newcommand{\hunum}{ h_{1} , \ldots ,h_{m} }
\newcommand{\iunuk}{ i_{1} , \ldots ,i_{k} }
\newcommand{\junul}{ j_{1} , \ldots ,j_{l} }
\newcommand{\lunum}{ l_{1} , \ldots ,l_{m} }
\newcommand{\llunum}{ \lambda_{1} , \ldots , \lambda_{m} }
\newcommand{\Lunuk}{ L_{1} , \ldots ,L_{k} }
\newcommand{\Punum}{ P_{1} , \ldots ,P_{m} }
\newcommand{\Punun}{ P_{1} , \ldots ,P_{n} }
\newcommand{\Qunum}{ Q_{1} , \ldots ,Q_{m} }
\newcommand{\Qunun}{ Q_{1} , \ldots ,Q_{n} }
\newcommand{\Runun}{ R_{1} , \ldots ,R_{n} }
\newcommand{\uunuk}{ u_{1} , \ldots ,u_{k} }

\newcommand{\aunun}{ a_{1}, \ldots ,a_{n} }
\newcommand{\aunudoi}{ a_{1}, a_{2} }
\newcommand{\bunudoi}{ b_{1}, b_{2} }
\newcommand{\ajunudoi}{ a_{j,1} , a_{j,2} }
\newcommand{\apunudoi}{ a_{1} p_{1} , p_{2} a_{2} }
\newcommand{\bunuunudoi}{ b_{1,1}, b_{1,2} }
\newcommand{\bhunudoi}{ b_{h,1} , b_{h,2} }
\newcommand{\bjunudoi}{ b_{j,1} , b_{j,2} }
\newcommand{\bkunudoi}{ b_{k,1} , b_{k,2} }
\newcommand{\bunun}{ b_{1}, \ldots ,b_{n} }

\newcommand{\funun}{ f_{1}, \ldots ,f_{n} }
\newcommand{\hlunum}{ (h_{1},l_{1}), \ldots ,(h_{m},l_{m}) }

\newcommand{\punudoi}{ p_{1} , p_{2} }
\newcommand{\punuunudoi}{ p_{1,1} , p_{1,2} }
\newcommand{\phunudoi}{ p_{h,1} , p_{h,2} }
\newcommand{\pjunudoi}{ p_{j,1} , p_{j,2} }
\newcommand{\pkunudoi}{ p_{k,1} , p_{k,2} }

\newcommand{\xunun}{ x_{1}, \ldots ,x_{n} }
\newcommand{\Xunun}{ X_{1}, \ldots ,X_{n} }
\newcommand{\zunun}{ z_{1}, \ldots ,z_{n} }
\newcommand{\zunudoi}{ z_{1} , z_{2} }
\newcommand{\zunuunudoi}{ z_{1,1} , z_{1,2} }
\newcommand{\zjunudoi}{ z_{j,1} , z_{j,2} }
\newcommand{\zkunudoi}{ z_{k,1} , z_{k,2} }

\newcommand{\ncps}{ ( {\cal A} , \varphi ) }
\newcommand{\noncom}{ \mbox{non-commutative probability space} }

\newcommand{\ncr}{ \mbox{non-crossing} }
\newcommand{\fps}{ \mbox{formal power series} }

\newcommand{\cprob}{ \mbox{$C^{*}$-probability space} }

\newcommand{\ncpol}{ \C \langle X_{1} , \ldots ,X_{n} \rangle } 
\newcommand{\coefis}{ \mbox{coef } ( \iunuk ) }
\newcommand{\seria}{ \sum_{k=1}^{\infty} \ 
     \sum_{i_{1}, \ldots ,i_{k} =1}^{n} \alpha_{(i_{1}, \ldots ,i_{k}) }
     z_{i_{1}}  \cdots  z_{i_{k}}   }

\newcommand{\ecpi}{ \stackrel{\pi}{\sim} }
\newcommand{\ecrho}{ \stackrel{\rho}{\sim} }
\newcommand{\eckinvrho}{ \stackrel{K^{-1} ( \rho )}{\sim} }
\newcommand{\ecsigma}{ \stackrel{\sigma}{\sim} }
\newcommand{\eccqsigma}{ \stackrel{\cq ( \sigma )}{\sim} }
\newcommand{\esigma}{ \stackrel{\crr ( \sigma )}{\sim} }
\newcommand{\eckpi}{ \stackrel{ K( \pi ) }{\sim} }
\newcommand{\egdef}{ \stackrel{def}{=} }
\newcommand{\ecdef}{ \stackrel{def}{ \Leftrightarrow } }

\newcommand{\nn}{ \{ 1, \ldots ,n \} }
\newcommand{\nnnn}{ \{ 1,2, \ldots ,2n \} }
\newcommand{\mm}{ \{ 1, \ldots ,m \} }
\newcommand{\mmmm}{ \{ 1,2, \ldots ,2m \} }
\newcommand{\kk}{ \{ 1, \ldots ,k \} }
\newcommand{\kd}{ K^{(2n)} }
\newcommand{\kdm}{ K^{(2m)} }
\newcommand{\cycle}{ ( 1 \rightarrow 2 \rightarrow \cdots \rightarrow
n \rightarrow 1 ) }
\newcommand{\ccycle}{ ( 1 \rightarrow 2 \rightarrow \cdots \rightarrow
2n \rightarrow 1 ) }

\newcommand{\freestar}{ \framebox[7pt]{$\star$} }

\newcommand{\cadj}{ c^{*} }
\newcommand{\uinv}{ u^{-1} }
\newcommand{\uadj}{ u^{*} }
\newcommand{\xadj}{ x^{*} }

\begin{document}

\title{\bf $R$-diagonal pairs - a common approach \\
to Haar unitaries and circular elements}
\author{Alexandru Nica  \thanks{Research done while this author was
on leave at the Fields Institute, Waterloo, and the Queen's University,
Kingston, holding a Fellowship of NSERC, Canada.}  \\
Department of Mathematics  \\
University of Michigan  \\
Ann Arbor, MI 48109-1003, USA  \\
(e-mail: andu@math.lsa.umich.edu) 
\and Roland Speicher  \thanks{Supported by a Heisenberg Fellowship
of the DFG.}  \\ 
Institut f\"{u}r Angewandte Mathematik  \\
Universit\"{a}t Heidelberg    \\
Im Neuenheimer Feld 294   \\
D-69120 Heidelberg, Germany   \\
(e-mail: roland.speicher@urz.uni-heidelberg.de) }
\date{November 1995}

\maketitle

\setlength{\baselineskip}{18pt}

{\large\bf Introduction} 

The Haar unitary and the circular element provide two of the most
frequently used $*$-distributions in the free probability theory
of Voiculescu, and in its remarkable applications to the study of
free products of von Neumann algebras (see e.g. [4,5,11,12,19,20]).
The starting point of the present paper is that if the $R$-transform
(i.e. the free analogue of the logarithm of the Fourier transform)
of these $*$-distributions is considered, then expressions of a
similar nature are obtained. The formulas are:
\[
\begin{array}{lc}
(I) &  
{ [ R( \mu_{u,  \uadj } ) ] ( z_{1} , z_{2} ) \ = \ 
\sum_{k=1}^{\infty} \frac{(-1)^{k+1} (2k-2)!}{(k-1)!k!}
(z_{1} z_{2} )^{k} + 
\sum_{k=1}^{\infty} \frac{(-1)^{k+1} (2k-2)!}{(k-1)!k!}
(z_{2} z_{1} )^{k} } 
\end{array}
\]
and
\[
\begin{array}{lccc}
(II) &  &  & 
{ [ R( \mu_{c, c^{*} } ) ] ( z_{1} , z_{2} ) \ = \ 
z_{1} z_{2} + z_{2} z_{1} , }
\end{array}
\]
where $u$ is a Haar unitary and $c$ is a circular element (in some
non-commutative probability spaces). Hence, a class of pairs of
elements which contains both $( u, \uadj )$ and $(c, \cadj )$ is
the following: for $x,y$ elements (random variables) in a
non-commutative probability space, the pair $(x,y)$ is said to be
{\em $R$-diagonal} if the $R$-transform $R( \mu_{x,y} )$ has the
form:
\[
\begin{array}{lccc}
(III) &  &  & 
{ [R( \mu_{x,y} )] ( \zunudoi ) \ = \ 
\sum_{k=1}^{\infty} \alpha_{k} (z_{1} z_{2})^{k} +
\sum_{k=1}^{\infty} \alpha_{k} (z_{2} z_{1})^{k} } 
\end{array}
\]
for some sequence $( \alpha_{k} )_{k=1}^{\infty}$ of complex
coefficients.

The main result of the paper (stated in Theorem 1.5, Corollary 1.8
below) is that the class of $R$-diagonal pairs has a remarkable
property of ``absorption'' under the operation of ``nested''
multiplication of free pairs. That is, for $a_{1}, a_{2}, p_{1},
p_{2}$ in a non-commutative probability space: if the sets
$\{ \aunudoi \}$ and $\{ \punudoi \}$ are free, and if the pair
$( \aunudoi )$ is $R$-diagonal, then so is $( a_{1}p_{1} ,
p_{2}a_{2} );$ and moreover, there exists a simple formula relating
the $\alpha$'s of Eqn.($III$) written for the two pairs
$( \aunudoi )$ and $( a_{1}p_{1} , p_{2}a_{2} ).$

As an immediate consequence, $R$-diagonal pairs exist in abundance,
even if we only want to consider pairs of the form $(x, x^{*}),$
and in the W$^{*}$-probabilistic context. We hope that, while not
as fundamental as $(u, \uadj )$ and $(c, \cadj )$ from Eqns.($I$),
($II$), other pairs in this class will also find their role in the 
theory, and in its applications.

A special interest is presented by the $R$-diagonal pairs of the
form $(up, (up)^{*} ),$ where $u$ is a Haar unitary and $p$ is
$*$-free from $u$ (in a $C^{*}$-probability space, say). A couple
of applications of the main result to this class is presented in
Sections 1.9, 1.10 below. A situation when several such free pairs  
$(up_{1} , ( up_{1} )^{*} ), \ldots , (up_{k} , ( up_{k} )^{*} )$
are considered at the same time is addressed in Theorem 1.13.

>From the technical point of view, our approach to the $R$-diagonal
pairs is based on the combinatorial description of the $R$-transform,
via the lattice $NC(n)$ of non-crossing partitions of $\nn ,$
$n \geq 1.$ More than once the proofs depend in an essential way
on considerations involving a certain operation $\freestar$ 
on formal power series, introduced in our previous paper 
\cite{NS2}; this operation represents in some sense the combinatorial 
facet of the $R$-transform approach to the multiplication of free
$n$-tuples of non-commutative random variables.

A detailed description of the results of the paper is made in the
next-coming Section 1. The rest of the paper is organized as follows.
In Section 2
we review some facts about non-crossing partitions, and in Section 3
we review the $R$-transform and the operation $\freestar$. Section 4
is devoted to pointing out a certain canonical bijection between
the set of intervals of $NC(n)$ and the set of 2-divisible
partitions in $NC(2n),$
which plays an important role in our considerations; we also note in 
Section 4 how, as a consequence of this combinatorial fact, one of the 
applications of $\freestar$ presented in \cite{NS2} can be improved.
The proofs of the results on $R$-diagonal pairs announced in Section 1 
are divided between the remaining Sections 5-8 of the paper. 

$\ $

$\ $

\setcounter{section}{1}

\setcounter{section}{1}
{\large\bf 1. Presentation of the results} 

$\ $

{\bf 1.1 Basic definitions}
In this section we briefly review some basic free probabilistic
terminology used throughout the paper (for a more detailed treatment,
we refer to the monograph \cite{VDN}).

\vspace{10pt}

{\bf The framework} We will call
{\em non-commutative probability space} a pair $\ncps$,
where $\A$ is a unital algebra (over $\C$), and 
$\varphi : \A \rightarrow \C$ is a linear functional 
normalized by $\varphi (1) =1.$ If we require in addition that 
$\A$ is a $C^{*}$-algebra, and $\varphi$ is positive, then
$\ncps$ is called a {\em $C^{*}$-probability space.} 

\vspace{10pt}

{\bf Freeness}
A family of unital subalgebras 
$\A_{1} , \ldots , \A_{n} \subseteq \A$ is said to be free in 
$\ncps$ if for every $k \geq 1,$ 
$1 \leq i_{1}, \ldots , i_{k} \leq n,$
$a_{1} \in \A_{i_{1}}, \ldots , a_{k} \in \A_{i_{k}}$, 
we have the implication:
\begin{equation}
\left\{  \begin{array}{c}
i_{1} \neq i_{2}, i_{2} \neq i_{3}, \ldots , i_{k-1} \neq i_{k}  \\
\varphi ( a_{1} ) =  \varphi ( a_{2} ) = \cdots = \varphi ( a_{k} ) =0
\end{array}  \right\}
\  \Rightarrow \ \varphi ( a_{1} a_{2} \cdots a_{k} ) =0.
\end{equation}
The notion of freeness in $\ncps$ extends to arbitrary subsets of 
$\A$, by putting $\X_{1} , \ldots , \X_{n} \subseteq \A$ to be free
if and only if the unital subalgebras generated by them are so. 
The freeness of a family of elements $\xunun \in \A$ is defined 
as the one of the family of subsets $\{ x_{1} \} , \ldots ,
\{ x_{n} \} \subseteq \A .$

If $\ncps$ is a $\cprob$, then the fact that $\xunun \in \A$ are 
$*$-free means by definition that the subsets 
$\{ x_{1} , x_{1}^{*} \} , \ldots , \{ x_{n} , x_{n}^{*} \}$ 
are free.

\vspace{10pt}

{\bf Joint distributions} The {\em joint distribution} 
of the family of elements 
$\aunun \in \A$, in the non-commutative probability space 
$\ncps$, is by definition the linear functional 
\newline
$\mu_{\aunun} : \ncpol \rightarrow \C$ given by:
\begin{equation}
\left\{  \begin{array}{l}
\mu_{\aunun} (1) \ = 1,  \\
\mu_{\aunun} ( \xiqxik ) \ = \ \varphi ( \aiqaik ) \ \ \mbox{ for }
k \geq 1 \mbox{ and } 1 \leq \iunuk \leq n, 
\end{array}  \right.
\end{equation}
where $\ncpol$ is the algebra of polynomials in $n$ non-commuting 
indeterminates $\Xunun$. 

If we only have one element $a = a_{1} \in \A$, then the functional
in (1.2) is just $\mu_{a} : \C [X] \rightarrow \C ,$
$\mu_{a} (f) = \varphi (f(a))$ for $f \in \C [X],$ and is called 
the {\em distribution} of $a$ in $\ncps .$

If $\ncps$ is a $\noncom$, then when speaking about the 
{\em $*$-distribution} of an element $x \in \A$, one usually
refers to the joint distribution $\mu_{x, x^{*}};$ another
(equivalent) approach goes by looking at the joint distribution
$\mu_{Re(x), \ Im(x)}$ of the real and imaginary parts
$Re(x) = (x+x^{*})/2,$ $Im(x) = (x-x^{*})/2i.$

The definitions of a Haar unitary and of a circular element
are made by prescribing (for the circular in an indirect way)
what is their $*$-distribution. We will consider the framework
of a $\cprob$ $\ncps$. Recall that:

- an element $u \in \A$ is called a {\em Haar unitary} in $\ncps$
if it is unitary, and if $\varphi (a^{n}) =0$ for every 
$n \in \Z \setminus \{ 0 \} ;$

- an element $a \in \A$ is called {\em semicircular} in $\ncps$ 
if it is selfadjoint and its distribution $\mu_{a}$ is
$\frac{1}{2 \pi} \sqrt{4-t^{2}} dt$ on [-2,2] (in other words, if
$\varphi (a^{n} ) = \frac{1}{2 \pi} \int_{-2}^{2} t^{n}
\sqrt{4-t^{2}} dt,$ $ n \geq 0);$  

- an element $c \in \A$ is called {\em circular} in $\ncps$
if it is of the form $(a+ib)/ \sqrt{2}$ with $a,b$ semicircular 
and free in $\ncps$. 

We take this occasion to review one more remarkable distribution:

- an element $a \in \A$ is called {\em quarter-circular} in $\ncps$
if it is positive and its distribution $\mu_{a}$ is
$\frac{1}{\pi} \sqrt{4-t^{2}} dt$ on [0,2] (i.e., if
$\varphi (a^{n} ) = \frac{1}{\pi} \int_{0}^{2} t^{n}
\sqrt{4-t^{2}} dt,$ $ n \geq 0).$  

\vspace{10pt}

{\bf The $R$-transform} Given a functional of the kind 
appearing in Eqn.(1.2) (i.e. $\mu : \ncpol \rightarrow \C$, 
such that $\mu (1) =1),$
its $R$-transform $R( \mu )$ is a certain formal power series 
in $n$ ``non-commuting complex variables'' $\zunun$:
\begin{equation}
[ R( \mu ) ] ( \zunun ) \ = \ \seria .
\end{equation}
The coefficients 
$( \alpha_{( \iunuk )} )_{k \geq 1, 1 \leq \iunuk \leq n}$
of $R( \mu )$ are also called the {\em free} (or
{\em non-crossing) cumulants} of $\mu .$
The precise definition of $R( \mu )$ (i.e. of how the free
cumulants are constructed from $\mu$) will be reviewed in Section 3
below.

If $\aunun$ are elements in the non-commutative probability space
$\ncps$, then the $R$-transform $R( \mu_{ \aunun } )$ contains the
information about the joint distribution $\mu_{\aunun}$ rearranged
in such a way that the freeness or non-freeness of  $\aunun$ becomes
transparent. More precisely, as proved in [14,8]
$\aunun$ are free in $\ncps$ if and only if the 
coefficient of $z_{i_{1}} z_{i_{2}} \cdots z_{i_{k}}$ in
$[ R( \mu_{\aunun} )] ( \zunun )$ vanishes whenever we don't have
$i_{1} = i_{2} = \cdots = i_{k} ;$ i.e., if and only if 
$R( \mu_{\aunun} )$ is of the form 
\begin{equation}
[ R( \mu_{\aunun} ) ] ( \zunun ) \ = \ 
f_{1} ( z_{1} ) + \cdots + f_{n} ( z_{n} ),
\end{equation}
for some formal power series of one variable $\funun .$

We will refer to the coefficient-vanishing condition presented 
in the preceding paragraph by saying that ``the series 
$R( \mu_{\aunun} )$ has no mixed coefficients''.  If this happens, then 
$\funun$ of (1.4) can only be the 1-dimensional $R$-transforms 
$R( \mu_{a_{1}} ) , \ldots , R( \mu_{a_{n}} )$, respectively.

$\ $

{\bf 1.2 Haar unitaries and circular elements}
If $F_{k}$ denotes the free group on generators $\gunuk$, then the
left-translation operators $\uunuk$ with $\gunuk$ on $l^{2} (F_{k})$
form a family of $*$-free Haar unitaries in $(L(F_{k} ) , \tau )$,
where $L(F_{k})$ is the von Neumann II$_{1}$ factor of $F_{k}$ and
$\tau$ is the unique normalized trace on $L(F_{k});$ $\uunuk$ is in
some sense ``the obvious system of generators'' for $L(F_{k}).$

In recent work of Voiculescu, Radulescu, Dykema (see e.g. 
[4,5,11,12,20]) it was shown that very
powerful results on $L(F_{k})$ can be obtained by using a different 
family of generators, consisting of free semicircular elements. The
Haar unitaries are also appearing in this picture, but in a more
subtle way, either via asymptotic models (for instance in \cite{V4},
Section 3), or via the theorem of Voiculescu \cite{V5} on the 
polar decomposition of the circular element. This theorem states 
that if $\A$ is a
von Neumann algebra, with $\varphi : \A \rightarrow \C$ a
faithful normal trace, and if $c$ is circular in $\ncps$, then 
by taking the polar decomposition $c=up$ of $c$ one gets that:
$u$ is a Haar unitary, $p$ is quarter-circular, and $u,p$ are
$*$-free. The original proof given by Voiculescu in \cite{V5}
for this fact depends on the asymptotic matrix model for free 
semicircular families developed in \cite{V4}. A direct, 
combinatorial proof was recently found by Banica \cite{Ba}.

The starting point of this work was the observation that the
$R$-transforms of the $*$-distributions of the Haar unitary
and of the circular element have similar forms. In fact, the
goal of the present paper is in some sense to understand
the relation between these two elements, from the point of 
view of the $R$-transform. Although this will not be our main 
concern, a new proof for the polar decomposition of the 
circular element will also follow (see the discussion in 
1.9, 1.10 below.)

If $c$ is a circular element in the
$\cprob$ $\ncps$, then from the fact that its real and imaginary
parts are multiples of free semicirculars, one gets immediately that:
\begin{equation}
[ R( \mu_{Re (c), \ Im (c)} ) ] ( z_{1} , z_{2} ) \ = \ 
\frac{1}{2} ( z_{1}^{2} + z_{2}^{2} );
\end{equation}
then by using a result from \cite{N} concerning linear changes of
coordinates, this implies:
\begin{equation}
[ R( \mu_{c, c^{*} } ) ] ( z_{1} , z_{2} ) \ = \ 
z_{1} z_{2} + z_{2} z_{1}.
\end{equation}
On the other hand, it was shown in [15, Section 3.4] that
for $u$ a Haar unitary in some $C^{*}$-probability space, one has: 
\begin{equation}
[ R( \mu_{u,  \uadj } ) ] ( z_{1} , z_{2} ) \ = \ 
\sum_{k=1}^{\infty} \frac{(-1)^{k+1} (2k-2)!}{(k-1)!k!}
(z_{1} z_{2} )^{k} + 
\sum_{k=1}^{\infty} \frac{(-1)^{k+1} (2k-2)!}{(k-1)!k!}
(z_{2} z_{1} )^{k} . 
\end{equation}

$\ $

Thus a class of pairs of elements which contains both 
$(c, \cadj )$ and $(u, \uadj )$ is given by the following

$\ $

{\bf 1.3 Definition:} Let $\ncps$ be a $\noncom$, and let
$\aunudoi$ be in $\A .$ We will say that $( \aunudoi )$
is an {\em $R$-diagonal pair} if 

(i) the coefficients of 
$\underbrace{z_{1} z_{2} \cdots z_{1} z_{2}}_{2n}$ and 
$\underbrace{z_{2} z_{1} \cdots z_{2} z_{1}}_{2n}$ in 
$[R( \mu_{ \aunudoi } )] ( \zunudoi )$ are equal, for
every $n \geq 1;$

(ii) every coefficient of $R( \mu_{\aunudoi} )$ not
of the form mentioned in (i) is equal to 0.

$\ $

In other words, the pair $( \aunudoi )$ is $R$-diagonal if 
and only if $R( \mu_{ \aunudoi } )$ has the form
\begin{equation}
[R( \mu_{\aunudoi} )] ( \zunudoi ) \ = \ 
\sum_{k=1}^{\infty} \alpha_{k} (z_{1} z_{2})^{k} +
\sum_{k=1}^{\infty} \alpha_{k} (z_{2} z_{1})^{k} 
\end{equation}
for some sequence $( \alpha_{k} )_{k=1}^{\infty}$. If this happens,
then the series of one variable $f(z) = \sum_{k=1}^{\infty}
\alpha_{k} z^{k}$ will be called the {\em determining series}
of the pair $( \aunudoi ).$

$\ $

{\bf 1.4 Remark} From (1.6) it follows that the determining 
series of $(c, \cadj )$ is just $f(z)=z.$ The determining series 
for $(u, \uadj )$, coming out from (1.7), will be denoted by $Moeb:$
\begin{equation}
Moeb (z) \ = \ \sum_{k=1}^{\infty} 
\frac{(-1)^{k+1} (2k-2)!}{(k-1)!k!} z^{k} ,
\end{equation}
and will be called the Moebius series (of one variable - compare
also to Eqn.(3.8) below).

$\ $

The main result of the paper can then be stated as follows.

$\ $

{\bf 1.5 Theorem} Let $\ncps$ be a $\noncom$, such that $\varphi$
is a trace (i.e. $\varphi (xy) = \varphi (yx),$ $x,y \in \A$),
and let $\aunudoi , \punudoi \in \A$ be such that $( \aunudoi )$ 
is an $R$-diagonal pair, and such that $\{ \punudoi \}$ is free from 
$\{ \aunudoi \} .$ Then $( a_{1} p_{1} , p_{2} a_{2} )$ is also 
an $R$-diagonal pair.

$\ $

Moreover, there exists a simple formula which connects the 
determining series of $( \aunudoi )$ and $( \apunudoi );$
this formula is presented in Corollary 1.8 below.

$\ $

{\bf 1.6 The operation $\freestar$} The formula announced in 
the previous phrase involves a certain binary operation 
$\freestar$ on the set of $\fps$ $\{ f \ |$
$f(z) = \sum_{k=1}^{\infty} \alpha_{k} z^{k} ; \ $
$\alpha_{1} , \alpha_{2} , \alpha_{3} , \ldots \in \C \} .$
One possible way of defining $\freestar$ is via the equation
\begin{equation}
R( \mu_{ab} ) \ = \ R( \mu_{a} ) \ \freestar \  R( \mu_{b} ),
\end{equation}
holding whenever $a$ is free from $b$ in some $\noncom$
$\ncps .$ (This definition makes sense because $R( \mu_{ab} )$
is completely determined by $R( \mu_{a} )$ and $R( \mu_{b} ),$ 
and because any two series $f,g$ of the considered type 
\footnote[1]{ In this paper the $R$-transform of an 1-dimensional
distribution is viewed as the particular case $n=1$ of the 
Eqn.(1.3); we warn the reader that this differs by a factor of
$z$ from the notation used in \cite{VDN}. }
can be realized as $R( \mu_{a} )$ and $R( \mu_{b} )$ with
$a,b$ free in some $\ncps .)$ The operation $\freestar$ also
has an alternative combinatorial definition, which will be 
reviewed in Section 3.3 below. The best point of view seems 
to be to consider both approaches to $\freestar$, 
and switch from one to the other as needed.
This operation appeared (under a different name) in
[14,9], in connection to the work of Voiculescu 
\cite{V3} on products of free elements. The name $\freestar$ 
was first used in \cite{NS2}, where the multivariable versions
of the operation were introduced and applied. We mention that
an important feature distinguishing the 1-dimensional instance
of $\freestar$ from the others is that in this (and only this)
case $\freestar$ is commutative.

$\ $

The formula announced immediately after Theorem 1.5 comes out 
in the following way.

$\ $

{\bf 1.7 Proposition} Let $\ncps$ be a $\noncom$, such that 
$\varphi$ is a trace, and let $( \aunudoi )$ be an $R$-diagonal 
pair in $\ncps$. If $f$ is the determining series of $( \aunudoi ),$ 
then: $f \ = \ R( \mu_{ a_{1} a_{2} } ) \ \freestar \ Moeb,$
where $Moeb$ is the Moebius series, as in Eqn.(1.9) (and of course,
$\mu_{ a_{1} a_{2} } : \C [X] \rightarrow \C$ denotes the 
distribution of the product $a_{1} \cdot a_{2} ).$

$\ $

{\bf 1.8 Corollary} In the context of Theorem 1.5, if $f$ and $g$
are the determining series of the $R$-diagonal pairs $( \aunudoi )$
and $( \apunudoi ),$ respectively, then we have the relation
\begin{equation}
g \ = \ f \ \freestar \ R( \mu_{ p_{1} p_{2} } ).
\end{equation}

$\ $

{\bf Proof} We can write that
\[
g \ = \ R( \mu_{ a_{1} p_{1} p_{2} a_{2} } )  \  \freestar \  Moeb 
\mbox{  (by Proposition 1.7) } 
\]
\[
 = \ R( \mu_{ a_{2} a_{1} p_{1} p_{2} } ) \ \freestar \  Moeb 
\mbox{  (because $\varphi$ is a trace) } \ 
\]
\[
= \ R( \mu_{ a_{2} a_{1} } ) \ \freestar \ R( \mu_{ p_{1} p_{2} } )
\ \freestar \ Moeb \mbox{  (by (1.10)) } 
\]
\[
= \ \left( R( \mu_{ a_{2} a_{1} } ) \ \freestar \ Moeb \right)
\ \freestar \  R( \mu_{ p_{1} p_{2} } )
\mbox{  (because $\freestar$ is commutative) } \ 
\]
\[
= \ f \ \freestar \ R( \mu_{ p_{1} p_{2} } )
\mbox{  (again by Proposition 1.7). {\bf QED} } 
\]

$\ $

{\bf 1.9 Application} Note that the Moebius series appears in two
different ways in the above discussion (in Remark 1.4 and 
Proposition 1.7, respectively). This has a consequence concerning
$R$-diagonal pairs of the form $(x, \xadj ),$ in the 
$C^{*}$-context. Let us denote by ${\cal R}_{c}$ the set of
formal power series of one variable which occur as $R( \mu ),$
with $\mu$ a probability measure on $\R$ having compact support
contained in $[0, \infty )$ (in connection to ${\cal R}_{c},$ see 
also \cite{V2}, Section 3). We have the 
following

\vspace{10pt}

{\bf Fact:} A formal power series $f$ of one variable can appear
as determining series for an $R$-diagonal pair $(x, \xadj )$ in
some $C^{*}$-probability space if and only if it is of the form
$f = g \freestar Moeb$, with $g \in {\cal R}_{c}.$ If this happens,
then $f$ can be in fact written as the determining series of
an $R$-diagonal pair $(up , (up)^{*} ),$ with $u$ Haar unitary,
$p$ positive, and such that $u$ is $*$-free from $p.$

\vspace{10pt}

{\bf Proof} Implication ``$\Rightarrow$'' follows from Proposition
1.7 ($f= R( \mu_{xx^{*}} ) \freestar Moeb$, and 
$R( \mu_{xx^{*}} ) \in {\cal R}_{c} ).$ Conversely, assume that
$f=g \freestar Moeb,$ with $g \in {\cal R}_{c}.$ We can always find
a $C^{*}$-probability space $\ncps$ and $u,p \in \A$, $*$-free,
such that $u$ is Haar unitary, $p$ is positive, and
$R( \mu_{p^{2}} )=g.$ (For instance we can take
$\A = L^{\infty} ( \mu ) \star L^{\infty} ( \T ),$ endowed with the
free product of $\mu$ with the Lebesgue measure on $\T$, where
$\mu$ is such that
$R( \mu ) =g.)$ The pair $(up , (up)^{*} )$ is $R$-diagonal by
Theorem 1.5, and has determining series $Moeb \freestar R( 
\mu_{p^{2}} ) = Moeb \freestar g =f,$ by Corollary 1.8.
{\bf QED}

\vspace{10pt}

We thus see that, from the point of view of the $*$-distribution,
any $R$-diagonal pair $(x, \xadj )$ in a $C^{*}$-probability space
can be replaced with one of the form $( up , (up)^{*} ).$ This
can be pushed to a ``polar decomposition result'', if we consider
the von Neumann algebra setting, with a normal faithful trace, and 
if we also assume
that Ker $x$ = $\{ 0 \} .$ Indeed, in such a situation we get that
the von Neumann subalgebras generated by $x$ and $up$ (in their
$W^{*}$-probability spaces) are canonically isomorphic, by an
isomorphism which sends $x$ into $up$ (this is the same type of
argument as, for instance, in \cite{V5}, Remark 1.10).
We obtain in this way a product decomposition $x =u'p',$
with $u'$ Haar unitary, $p'$ positive, and $u'$ $*$-free from $p',$
and the uniqueness of the polar decomposition shows that $u'p'$ is
necessarily the polar decomposition of $x.$ (The needed fact that
Ker $p'$ = $\{ 0 \}$ is obtained by verifying that the distribution of
$p'^{2}$ has no atom at 0.)

$\ $

{\bf 1.10 Application} Let $\ncps$ be a $\cprob$, with $\varphi$ a
trace, and let $u,p \in \A$ be such that $u$ is a Haar unitary,
$*$-free from $p.$ We look for necessary and sufficient
conditions for the real and imaginary parts of $up$ to be free.
Without loss of generality, we can assume that $p$ is normalized
in such a way that $\varphi (pp^{*} ) =1.$

By using 1.5, 1.8 and also Eqn.(1.7) of 1.2 we infer that 
$(up, p^{*} u^{*} )$ is an $R$-diagonal pair, with determining
series $g$ given by:
\begin{equation}
g \ = \ Moeb \ \freestar \ R( \mu_{ pp^{*} } ) .
\end{equation}
We write explicitly $g(z) = \sum_{k=1}^{\infty} \beta_{k} z^{k} ;$
the assumption that $\varphi (pp^{*} ) =1$ plugged into (1.12)
implies that $\beta_{1} =1.$ Remembering how the determining series
was defined in 1.3, we have that
\[
[ R( \mu_{up , p^{*} u^{*} } ) ] ( \zunudoi ) \ = \ 
\sum_{k=1}^{\infty} \beta_{k} ( z_{1} z_{2} )^{k}  \ + \
\sum_{k=1}^{\infty} \beta_{k} ( z_{2} z_{1} )^{k}  ;
\]
then by doing a linear change of coordinates (as in \cite{N},
Section 5) we get
\begin{equation}
[ R( \mu_{Re(up) , Im(up) } ) ] ( \zunudoi ) \ = \ 
\sum_{k=1}^{\infty} \frac{\beta_{k}}{4^{k}}
( ( z_{1} + i z_{2} )( z_{1} - i z_{2}) )^{k}  \ + \
\sum_{k=1}^{\infty} \frac{\beta_{k}}{4^{k}}
( ( z_{1} - i z_{2} )( z_{1} + i z_{2}) )^{k}  .
\end{equation}

By the result stated in (1.4), $Re(up)$ and $Im(up)$ are free if and only
if the series in (1.13) has no mixed coefficients; but a direct analysis 
of the right-hand side of (1.13) shows that this can happen if and 
only if $\beta_{2} = \beta_{3} = \cdots =0.$ Hence $Re(up)$ and
$Im(up)$ are free if and only if the series $g$ of (1.12) is just
$g(z) =z.$ Finally, the equation
$[ Moeb \ \freestar \ R( \mu_{pp^{*}} ) ] (z) =z$ is easily solved
``in the unknown'' $R( \mu_{ pp^{*} } ),$ and is found to be
equivalent to $[R( \mu_{ pp^{*} } )] (z) = z/(1-z).$
We thus obtain the following

\vspace{10pt}

{\bf Fact:} With $u$ and $p$ as in the first paragraph of 1.10,
the necessary and sufficient condition for the real and imaginary 
part of $up$ to be free is
\begin{equation}
[R( \mu_{ pp^{*} } )] (z) = z/(1-z).
\end{equation}

\vspace{10pt}

Note that if $\beta_{1} =1$ and $\beta_{2} = \beta_{3} = \cdots
= 0,$ then the right-hand side of (1.13) is just 
$(z_{1}^{2} + z_{2}^{2} )/2 ,$ and by comparing (1.13) against 
(1.5) we find that $up$ is circular. Hence $up$ is circular
whenever (1.14) holds. 

We also note that (1.14) does hold if $p=p^{*} =$ 
quarter-circular (the square of the quarter-circular is the same 
thing as the square of the semicircular, and the $R$-transform of 
the latter square is well-known - see e.g. \cite{NS2}, Lemma 4.2 
or Lemma 1.1 in the Appendix). Here again
the uniqueness of the polar decomposition yields from this point 
a proof for the polar decomposition of the circular element. 

$\ $

{\bf 1.11 Case of several pairs}
Another aspect which can be studied in the context of 1.10 is:
what happens if instead of looking just at $up,$ we look at 
a family $up_{1} , up_{2}, \ldots , up_{k},$ where
$u, p_{1}, p_{2} , \ldots , p_{k}$ are $*$-free? It is shown
by Banica in \cite{Ba} that $up_{1} , up_{2}, \ldots , up_{k}$ 
are also $*$-free if the following happens: every $p_{j}$
$(1 \leq j \leq k)$ is in some sense $*$-modeled by an operator
of the form $S^{*} + S^{m_{j}},$ $m_{j} \geq 1,$ where $S$ is
the unilateral shift on $l^{2} ( \N ).$ We have found a
condition expressed in terms of ``moments of pairs'', which 
still implies the $*$-freeness of $up_{1} , up_{2}, \ldots , up_{k}.$ 
This condition, which is satisfied by the pair
$( S^{*} + S^{m} , S + {S^{m}}^{*} )$ for every $m \geq 1 ,$
is described as follows.

$\ $

{\bf 1.12 Definition} Let $\ncps$ be a $\noncom$, and let
$\aunudoi$ be in $\A .$ We will say that $( \aunudoi )$ is a
{\em diagonally balanced pair} if
\begin{equation}
\varphi ( \underbrace{a_{1}a_{2} \cdots a_{1}a_{2}a_{1} }_{2n+1} ) =
\varphi ( \underbrace{a_{2}a_{1} \cdots a_{2}a_{1}a_{2} }_{2n+1} ) =
0, \ \mbox{ for every } n \geq 0.
\end{equation}

$\ $

Then we have:

$\ $

{\bf 1.13 Theorem} Let $\ncps$ be a $\noncom$, with $\varphi$
a trace, and let $u, p_{1,1}, p_{1,2} , \ldots , p_{k,1},
p_{k,2}$ be in $\A$ such that:

(i) $u$ is invertible and $\varphi (u^{n}) =0$ for every
$n \in \Z \setminus \{ 0 \};$

(ii) the pairs $( p_{1,1} , p_{1,2} ), \ldots , 
( p_{k,1} , p_{k,2} )$ are diagonally balanced;

(iii) the sets $\{ u, \uinv \} ,
\{ p_{1,1}, p_{1,2} \} , \ldots , \{ p_{k,1} , p_{k,2} \}$
are free.
\newline
Then the sets
$\{ up_{1,1}, p_{1,2} \uinv \} , \ldots , 
\{ up_{k,1}, p_{k,2} \uinv \}$ are also free.

$\ $

The techniques used for proving the results announced in 1.5,
1.7, 1.13 above are based on the combinatorial approach to the
$R$-transform, via the lattice $NC(n)$ of non-crossing partitions
of $\nn$, $n \geq 1.$ In some instances, we are able to use an 
elegant idea of Biane \cite{Bi} which reduces assertions on
non-crossing partitions to calculations in the group algebra of
the symmetric group (but there are also situations when this
mechanism is apparently not applying, and we have to use
``geometric'' arguments). It seems that from the combinatorial
point of view, a certain canonical bijection between the
set of intervals of $NC(n)$ and the set of 2-divisible
partitions in $NC(2n)$ has a significant role in the considerations.
We have incidentally noticed that, as a consequence of this
combinatorial fact, we can improve one of the applications of
the operation $\freestar$ that were presented in \cite{NS2}.
Namely, we have:

$\ $

{\bf 1.14 Theorem} Let $\ncps$ be a $\noncom$, with $\varphi$
a trace, and let $\aunun , b' , b'' \in \A$ be such that the pair 
$(b',b'')$ is diagonally balanced, and such that $\{ b',b'' \}$
is free from $\{ \aunun \} .$ 
Then $\{ b'a_{1}b'' , \ldots , b'a_{n}b'' \}$ is free from 
$\{ \aunun \} .$

$\ $

The particular case of 1.14 when $b'$ = $b''$ is a 
semicircular element is stated in Application 1.10 of \cite{NS2}.
The Theorem 1.14 sounds quite ``elementary'', and it is not
impossible that it also has a simple direct proof, using only 
the definition of freeness (we weren't able to find one, though).

$\ $

$\ $

$\ $

\setcounter{section}{2}
{\large\bf 2. Preliminaries on non-crossing partitions} 

$\ $

{\bf 2.1 Definition of NC(n)} If $\pi = \{ B_{1} , \ldots , B_{r} \}$
\setcounter{equation}{0}
is a partition of $\nn$ (i.e. $B_{1} , \ldots , B_{r}$ are pairwisely
disjoint, non-void sets, such that $B_{1} \cup \cdots \cup B_{r}$ =
$\nn$), then the equivalence relation on $\nn$ with equivalence classes
$B_{1} , \ldots , B_{r}$ will be denoted by $\ecpi$; the sets  
$B_{1} , \ldots , B_{r}$ will be also referred to as the {\em blocks}
of $\pi$. The number of elements in the block $B_k$, $1 \leq k \leq r,$
will be denoted by $|B_{k}|$.

A partition $\pi$ of $\nn$ is called {\em non-crossing} 
if for every $1 \leq i < j < i' < j' \leq n$ such that $i \ecpi i'$
and $j \ecpi j'$, it necessarily follows that
$i \ecpi j \ecpi i' \ecpi j'$. The set of all non-crossing 
partitions of $\nn$ will be denoted by $NC(n).$ On $NC(n)$ we will
consider the {\em refinement order,} defined by
$\pi \leq \rho \ \ecdef$ each block of $\rho$ is a union of blocks 
of $\pi$ (in other words, $\pi \leq \rho$ means that the implication 
$i \ecpi j \Rightarrow i \ecrho j$ holds, for $1 \leq i,j \leq n).$
The partially ordered set $NC(n)$ was introduced by G. Kreweras in
\cite{K}, and its combinatorics has been studied by several authors
(see e.g. \cite{SU}, and the list of references there). 

$\ $

{\bf 2.2 The circular picture} of a partition  $\pi = \{ B_{1},
\ldots , B_{r} \}$ of $\nn$ is obtained by drawing $n$ equidistant
and clockwisely ordered points $\Punun$ on a circle, and then by
drawing for each block $B_{k}$ of $\pi$ the inscribed convex polygon
with vertices $\{ P_{i} \ | \ i \in B_{k} \}$ (this polygon may
of course be reduced to a point or a line segment). It is immediately
verified that $\pi$ is non-crossing if and only if the $r$ convex 
polygons obtained in this way are disjoint.

We take the occasion to mention that when $B$ is a block of the 
partition $\pi \in NC(n),$ we will use for $i<j$ in $B$ the 
expression ``$i$ and $j$ are consecutive in $B$'' to mean that
either $B \cap \{ i+1, \ldots , j-1 \} = \emptyset$ or $i$ =
min $B,$ $j$ = max $B.$ On the circular picture, the fact that $i,j$
are consecutive in $B$ means that $P_{i} P_{j}$ is an edge (rather
than a diagonal) of the inscribed polygon with vertices 
$\{ P_{h} \ | \ h \in B \} .$

$\ $

{\bf 2.3 The  Kreweras complementation map} is a remarkable order
anti-isomorphism $K:NC(n) \rightarrow NC(n)$, introduced in
$\cite{K}$, Section 3, and described as follows. 

Let $\pi$ be in $NC(n),$ and consider the circular picture of 
$\pi ,$ involving the points $P_{1} , \ldots , P_{n},$ as in 2.2.
Denote the midpoints of the arcs of circle $P_1 P_2 , \ldots , 
P_{n-1} P_n , \ P_n P_1$ by $Q_{1}, \ldots , Q_{n-1},$
\newline
$ Q_{n}$, respectively. Then the ``complementary'' partition $K( \pi ) 
\in NC(n)$ is given, in terms of the corresponding equivalence 
relation $\eckpi$ on $\nn$, by putting for $1 \leq i,j \leq n:$
\begin{equation}
i \eckpi j \ \ecdef \ \left\{  \begin{array}{l}
{ \mbox{ there are no $1 \leq h,k \leq n$ such that 
$h \ecpi k$ and such that}  }  \\
{ \mbox{ the line segments $Q_{i} Q_{j}$ and $P_{h} P_{k}$ have 
non-void intersection.} } 
\end{array}   \right.
\end{equation}
It is easily verified that $\eckpi$ of (2.1) is indeed an equivalence
relation on $\nn$, and that the partition $K( \pi )$ corresponding to 
it is non-crossing (for the latter thing, we just have to look at the
circular picture of $K( \pi ),$ with respect to the points $\Qunun ).$

In order to verify that $K: NC(n) \rightarrow NC(n)$ defined in this
way really is an order anti-isomorphism, one can first remark that
$K^{2} ( \pi )$ is (for every $\pi \in NC(n)$) a
rotation of $\pi$ with $360^{o}/n$; this shows in particular that 
$K$ is a bijection. The implication 
$\pi \leq \rho \Rightarrow K( \pi ) \geq K( \rho )$  is a direct
consequence of (2.1), and the converse must also hold, since 
$K^{2}$ is an order-preserving isomorphism of $NC(n).$

As a concrete example, the circular-picture verification for 
$K( \{ \{ 1,4,5 \} , \{ 2,3 \} ,$
$\{ 6,8 \} ,$
\newline
$\{ 7 \} \} )$ =
$\{ \{ 1,3 \} , \{ 2 \} , \{ 4 \} , \{ 5,8 \} , \{ 6,7 \} \}
\in NC(8)$ is shown in Figure 1.
\newpage
$\ $

\vspace{5cm}

\[
\mbox{\bf Figure 1.}
\]

By examining the common circular picture for $\pi$ and $K( \pi )$,
it becomes quite obvious that $K( \pi )$ could be alternatively
defined as the biggest $\rho$ in $(NC(n) , \leq )$ with the 
property that the partition of $\nnnn$ obtained by interlacing 
$\pi$ and $\rho$ is still non-crossing. This fact is formally
recorded in the next proposition.

$\ $

{\bf 2.4 Proposition} Let $\pi$ and $\rho$ be in $NC(n)$. Denote by
$\pi '$ and $\rho '$ the partitions of $\{ 2,4, \ldots ,$
$2n \}$ and $\{ 1,3, \ldots , 2n-1 \}$, respectively, which get 
identified to $\pi$ and $\rho$ via the order-preserving bijections 
$ \nn \rightarrow \{ 2,4, \ldots ,2n \}$ and  
$ \nn \rightarrow \{ 1,3, \ldots ,2n-1 \}$. Denote by $\sigma$ the 
partition of $\{ 1,2, \ldots , 2n \}$ formed by putting $\pi '$ and 
$\rho '$ together. Then $\sigma$ is $\ncr$ if and only if 
$\pi \leq K( \rho ).$

$\ $

{\bf 2.5 The relative Kreweras complement} Given $\rho \in NC(n),$
one can define a relativized version of the Kreweras complementation
map, $K_{\rho},$ which is an order anti-isomorphism of
$\{ \pi \in NC(n) \ | \ \pi \leq \rho \} .$
If we write explicitly $\rho = \{ B_{1}, \ldots , B_{r} \} ,$ then an
arbitrary element of
$\{ \pi \in NC(n) \ | \ \pi \leq \rho \}$ can be written as
$\{ A_{1,1}, \ldots , A_{1,s_{1}}, \ldots , A_{r,1}, \ldots , 
A_{r,s_{r}} \} ,$ where $A_{k,1} \cup \cdots \cup A_{k,s_{k}} =
B_{k},$ $1 \leq k \leq r.$ The relative Kreweras complement 
$K_{\rho} ( \pi )$ is obtained by looking for every $1 \leq k \leq r$
at the non-crossing partition 
$\{ A_{k,1} , \ldots , A_{k,s_{k}} \}$ of $B_{k},$ and by
taking its Kreweras complement, call it $\theta_{k} ,$ in the sense
of Section 2.3 above (of course, in order to do this, we need to
identify canonically $B_k$ with $\{ 1, \ldots , |B_{k}| \} )$.
Then $K_{\rho} ( \pi )$ is obtained by putting together the
partitions $\theta_{1}$ of $B_{1},$ $\ldots$ , $\theta_{r}$ of
$B_{r}$ (see also \cite{NS2}, Sections 2.4, 2.5).

Note that if $\rho = \{ \nn \}$ is the maximal element of
$(NC(n), \leq )$, then $K_{\rho}$ coincides with
$K:NC(n) \rightarrow NC(n)$ discussed in 2.3.

$\ $

{\bf 2.6 Relation with permutations} A useful way of
``encoding'' non-crossing partitions by permutations was
introduced by Ph. Biane in \cite{Bi}. Let ${\cal S}_n$
denote the group of all permutations of $\nn$. For 
$B = \{ i_1 < i_2 < \cdots < i_m \} \subseteq \nn$
we denote by $\gamma_{B} \in {\cal S}_n$ the cycle given by
\begin{equation}
\left\{  \begin{array}{l}
{ \gamma_{B} ( i_1 ) = i_2 , \ldots , \gamma_{B} ( i_{m-1} ) = i_m , 
\gamma_{B} ( i_{m} ) = i_1 ,}   \\
{ \gamma_{B} (j) = j \mbox{ for } j \in \nn \setminus B }
\end{array}  \right.
\end{equation}
(if $|B| =1,$ we take $\gamma_{B}$ to be the unit of ${\cal S}_n$).
Then for every $\pi \in NC(n),$ the permutation associated to it is
\begin{equation}
Perm ( \pi ) \  = \ \prod_{B \ block \ of \ \pi}  \gamma_{B} 
\end{equation}
(the cycles $( \gamma_{B} \ ; \ B$ block of $\pi$) commute, so
the order of the factors in the product (2.3) does not matter).

It was shown in \cite{Bi} how the (obviously injective) map
$Perm : NC(n) \rightarrow {\cal S}_{n}$ can be used for an
elegant analysis of the skew-automorphisms (i.e. automorphisms
or anti-automorphisms) of $(NC(n) , \leq ).$ We will only need
here the ``$Perm$'' characterization of the Kreweras complementation
map, which goes as follows: if $\gamma \in {\cal S}_{n}$ denotes 
the cycle $\cycle$, then  
\begin{equation}
Perm (K( \pi )) \ = \ Perm ( \pi )^{-1} \ \gamma ,
\mbox{  for every } \pi \in NC(n). 
\end{equation} 
This characterization can be extended without difficulty to the
situation of the relative Kreweras complement, we have:
\begin{equation}
Perm ( K_{\rho} ( \pi )) \ = \ Perm ( \pi )^{-1} \ Perm ( \rho ),
\end{equation} 
for every $\pi , \rho \in NC(n)$ such that $\pi \leq \rho$
(see \cite{NS2}, Section 2.5).

$\ $

$\ $

$\ $

\setcounter{section}{3}
{\large\bf 3. Review of the operation} $\freestar$ 
{\large\bf and of the} $R${\large\bf -transform} 

$\ $

We will follow the presentation of \cite{NS2}, Section 3. The 
\setcounter{equation}{0}
approach to the $R$-transform taken here is based on elements of
Moebius inversion theory for non-crossing partitions, on the lines
of \cite{S1}. We mention that the connection between this approach 
and the original one of Voiculescu in \cite{V2} is made via an
alternative description of the $n$-dimensional $R$-transform, which
goes by ``modeling on the full Fock space over $\C^{n}$'' (see
\cite{N}).

$\ $

{\bf 3.1 Notation} Let $n$ be a positive integer. We denote by 
$\Theta_{n}$ the set of formal power series without 
constant coefficient in $n$ non-commuting variables 
$z_1 , \ldots , z_n .$ An element of $\Theta_{n}$ is thus a
series of the form
\begin{equation}
f( \zunun ) \ = \ \seria ,
\end{equation} 
where $( \alpha_{( \iunuk )} \ ; \ k \geq 1, \ 1 \leq \iunuk 
\leq n)$ is a family of complex coefficients.

$\ $

{\bf 3.2 Notations for coefficients} The following conventions 
for denoting coefficients of formal power series will be used
throughout the whole paper.

$1^{o}$ For $f \in \Theta_{n}$ and
$k \geq 1$, $1 \leq \iunuk \leq n,$ we will denote
\begin{equation}
[ \coefis ] (f) \ \egdef \  \mbox{the coefficient of $\ziqzik$ in $f$} .
\end{equation}

$2^{o}$ Restrictions of $k$-tuples: let $k \geq 1$ and 
$1 \leq \iunuk \leq n$ be integers, and let $B = \{ h_{1} <
h_{2} < \cdots < h_{r} \}$ be a non-void subset of $\kk .$
Then by ``$( \iunuk ) | B$'' we will understand the $r$-tuple 
$(i_{h_{1}} , i_{h_{2}} ,  \ldots , i_{h_{r}} ).$ An expression 
like
\begin{equation}
[ \mbox{coef } ( \iunuk )|B ] (f)
\end{equation}
for $f \in \Theta_{n},$ $k \geq 1, \ 1 \leq \iunuk \leq n$ and
$\emptyset \neq B \subseteq \kk$, will hence mean that the 
convention of notation (3.2) is applied to the 
$|B|$-tuple $( \iunuk )|B.$

$3^{o}$ Given $f \in \Theta_{n},$ $k \geq 1, \ 1 \leq \iunuk \leq n$
integers, and $\pi$ a non-crossing partition of $\kk ,$ we will denote
\begin{equation}
[ \coefis ; \pi ] (f) \ \egdef \ \prod_{B \ block \ of \ \pi}
[ \coefis |B ] (f) .
\end{equation}
Thus if $n,k, \ 1 \leq \iunuk \leq n$ and $\pi$ are fixed, then
$[ \coefis ; \pi ] : \Theta_{n} \rightarrow \C$ is a functional,
generally non-linear (it is linear if and only if 
$\pi$ is the partition into only one block, $\{ \kk \} ,$ in which
case $[ \coefis ; \pi ] = [ \coefis ]$ of (3.2)). 

$4^{o}$ The 1-dimensional case: If $n=1,$ then the convention of
notation in (3.2) can (and will) be abridged from 
$[ \mbox{coef } ( \underbrace{1, \ldots ,1}_{k} ) ] (f)$ to just
$[ \mbox{coef } (k)] (f).$ A similar abbreviation will be used for 
the convention in (3.4); note that the restriction of $k$-tuples 
gets in this case the form ``$(k)|B = ( |B| )$'', for 
$\emptyset \neq B \subseteq \kk ,$ hence (3.4) becomes:
\begin{equation}
[ \mbox{coef } (k); \pi ] (f) \ = \ 
\prod_{B \ block \ of \ \pi} [ \mbox{coef } (|B|) ] (f),
\end{equation}
for $f \in \Theta_{1}, \ k \geq 1$ and $\pi \in NC(k).$

$\ $

{\bf 3.3 The operation \freestar} Let $n$ be a positive integer.
We denote by $\freestar \ (= \freestar_{n} )$ the binary 
operation on the set $\Theta_{n}$ of 3.1, determined by the formula
\begin{equation}
[ \coefis ] (f \freestar g) \ = \  \sum_{\pi \in NC(k)}
[ \coefis ; \pi ] (f) \cdot [ \coefis ; K( \pi ) ] (g) ,
\end{equation}
holding for every $f,g \in \Theta_{n}, \ k \geq 1, \ 1 \leq
\iunuk \leq n,$ and where $K: NC(k) \rightarrow NC(k)$ is the 
Kreweras complementation map reviewed in Section 2.3.

The operation $\freestar$ on $\Theta_{n}$ is associative
(see \cite{NS2}, Proposition 3.5) and is commutative when (and 
only when) $n=1$ (see e.g. \cite{NS1}, Proposition 1.4.2).
The unit for $\freestar$ is the series which takes the sum of
the variables, $Sum ( \zunun ) \egdef z_{1} + \cdots + z_{n}.$

An important role in the considerations related to $\freestar$
is played by the series
\begin{equation}
Zeta ( \zunun ) \ = \ \sum_{k=1}^{\infty} \ 
\sum_{\iunuk =1}^{n}  \ziqzik
\end{equation}
and
\begin{equation}
Moeb ( \zunun ) \ = \ \sum_{k=1}^{\infty} \ \sum_{\iunuk =1}^{n}  
(-1)^{k+1}  \frac{ (2k-2)! }{ (k-1)! k! } \ \ziqzik , 
\end{equation}
which are called the ($n$-variable) Zeta and Moebius series, 
respectively. (The names are coming from the combinatorial 
interpretation of these series. The relation with the Moebius inversion
theory in a poset, as developed in \cite{DRS}, is particularly clear
in the case when $n=1$ - see \cite{NS1}.) $Zeta$ and $Moeb$ are inverse
to each other with respect to $\freestar ,$ i.e. $Zeta \freestar Moeb$ =
$Sum$ = $Moeb \freestar Zeta.$ Moreover, they are central
with respect to $\freestar$ (i.e. $Zeta \freestar f = f \freestar Zeta$
for every $f \in \Theta_{n},$ and similarly for $Moeb).$

$\ $

Although chronologically the $R$-transform preceded the 
$\freestar$-operation, it is convenient to define it here in the 
following way.

$\ $

{\bf 3.4 Definition} Let $n$ be a positive integer. To every
functional $\mu : \ncpol \rightarrow \C ,$ normalized by
$\mu (1)=1$ (where $\ncpol$ is the algebra of polynomials in
$n$ non-commuting indeterminates, as in (1.2)), we attach two
formal power series $M( \mu ), R( \mu ) \in \Theta_{n},$ by
the equations:
\begin{equation}
[M( \mu )] ( \zunun ) \ = \ \sum_{k=1}^{\infty} \
\sum_{\iunuk =1}^{n} \mu ( \xiqxik ) \ziqzik
\end{equation}
and
\begin{equation}
R( \mu ) \ = \ M( \mu ) \ \freestar \ Moeb. 
\end{equation}
The series $R( \mu )$ is called the $R$-transform of $\mu .$

$\ $

{\bf 3.5 The moment-cumulant formula} Since $Moeb$ is the inverse
of $Zeta$ under $\freestar ,$ the Equation (3.10) is clearly
equivalent to
\begin{equation}
M( \mu ) \ = \ R( \mu ) \ \freestar \ Zeta;
\end{equation}
(3.11) is in some sense ``the formula for the $R^{-1}$-transform'', 
since the transition from $M( \mu )$ back to $\mu$ is trivial.
The formula (3.11) is more frequently used than (3.10), for the
reason that the coefficients of $Zeta$ are easier to handle than 
those of $Moeb.$ The equation obtained by plugging (3.6) into (3.11) 
and taking into account that $[ \mbox{coef } ( \iunuk ); \pi ] (Zeta)$
is always equal to 1 was first observed in \cite{S1}, and looks
like this:
\begin{equation}
[ \mbox{coef } ( \iunuk ) ] (M( \mu )) \ = \ 
\sum_{\pi \in NC(k)} 
[ \mbox{coef } ( \iunuk ); \pi ]  (R( \mu )),
\end{equation}
for $k \geq 1$ and $1 \leq \iunuk \leq n.$ We will call (3.12)
the {\em moment-cumulant formula,} because it connects the 
coefficient $[ \mbox{coef } ( \iunuk ) ] (M( \mu )) $
- which is just $\mu ( \xiqxik ),$ a ``moment'' of $\mu ,$
with the coefficients of $R( \mu )$ - which were called in 
\cite{S1} the ``free (or non-crossing) cumulants'' of $\mu .$

$\ $

The multi-variable $R$-transform has the fundamental property 
that it stores the information about a joint distribution (in the
sense of Eqn.(1.2)) in such a way that the freeness or
non-freeness of the non-commutative random variables involved 
becomes very transparent. More precisely, we have the following

$\ $

{\bf 3.6 Theorem [14,8]:} The families of elements
$\{ a_{1,1}, \ldots ,a_{1,m_{1}} \} ,$ $\ldots$ ,
$\{ a_{n,1}, \ldots ,a_{n,m_{n}} \}$ are free in the non-commutative
probability space $\ncps$ if and only if the coefficient of
$z_{i_{1},j_{1}} z_{i_{2},j_{2}} \cdots z_{i_{k},j_{k}}$ in
$[ R( \mu_{a_{1,1}, \ldots ,a_{1,m_{1}} , \ldots , a_{n,1}, 
\ldots ,a_{n,m_{n}}} )]$
$( z_{1,1}, \ldots ,z_{1,m_{1}} , \ldots , z_{n,1}, 
\ldots ,z_{n,m_{n}} )$
vanishes whenever we don't have $i_{1} = i_{2} = \cdots = i_{k};$
i.e. if and only if 
$R( \mu_{a_{1,1}, \ldots ,a_{1,m_{1}} , \ldots , a_{n,1}, 
\ldots ,a_{n,m_{n}}} )$ is of the form
\begin{equation}
f_{1} ( z_{1,1}, \ldots ,z_{1,m_{1}} ) + \cdots +  
f_{n} ( z_{n,1}, \ldots ,z_{n,m_{n}} )  
\end{equation}
for some formal power series $\funun .$

$\ $

We will refer to the coefficient-vanishing condition in Theorem 3.6 by
saying that ``the $R$-series has no mixed coefficients''. If
this happens, then $\funun$ of (3.13) can only be the 
$R$-transforms $R( \mu_{a_{1,1}, \ldots ,a_{1,m_{1}}} )$, 
$\ldots$ , $R( \mu_{a_{n,1}, \ldots ,a_{n,m_{n}}} ),$ respectively. 

$\ $

The Equation (3.6) of 3.3 will be called in the sequel ``the
combinatorial definition of $\freestar$''. In an approach where
the $R$-transform is defined first, another possible definition
of $\freestar$ would be provided by the following fact (which
comes here as a theorem).

$\ $

{\bf 3.7 Theorem (\cite{NS2})} Let $\ncps$ be a $\noncom$, and
let $a_{1} , \ldots , $
\newline
$a_{n}, b_{1} , \ldots b_{n} \in \A$ be 
such that $\{ \aunun \}$ is free from $\{ \bunun \} .$ Then
\begin{equation}
R( \mu_{a_{1}b_{1} , \ldots , a_{n}b_{n}} ) \ = \ 
R( \mu_{\aunun} ) \  \freestar \ R( \mu_{\bunun} ).
\end{equation}

$\ $

We mention that it is often handier to use the Equation (3.14)
after performing a $\freestar$-operation with $Zeta$ on the right
on both sides, and invoking (3.11); this leads to the equivalent
equations:
\begin{equation}
M( \mu_{a_{1}b_{1} , \ldots , a_{n}b_{n}} ) \ = \ 
R( \mu_{\aunun} ) \  \freestar \ M( \mu_{\bunun} )
\end{equation}
or
\begin{equation}
M( \mu_{a_{1}b_{1} , \ldots , a_{n}b_{n}} ) \ = \ 
M( \mu_{\aunun} ) \  \freestar \ R( \mu_{\bunun} )
\end{equation}
(obtained by attaching the $Zeta$ on the right-hand side to
$R( \mu_{\bunun} )$ and to $R( \mu_{\aunun} ),$ respectively).

$\ $

Finally, since in this paper we are dealing only with tracial
non-commutative probability spaces, it is useful to recall the
following simple fact concerning cyclic permutations in the
$R$-transform.

$\ $

{\bf 3.8 Proposition} Let $n$ be a positive integer and let
$\mu : \ncpol \rightarrow \C$ be a linear functional, such that
$\mu (1)=1.$ Assume that $\mu$ is a trace (i.e. $\mu (p'p'') =
\mu (p''p')$ for every two polynomials $p',p'' \in \ncpol ;$
this is for instance the case whenever $\mu = \mu_{\aunun}$ for
some elements $\aunun$ in some $\ncps$, with $\varphi$ a trace).
Then the coefficients of $R( \mu )$ are invariant under cyclic
permutations, i.e.
\begin{equation}
[ \mbox{coef } ( \iunuk , \junul ) ] (R( \mu )) \ = \ 
[ \mbox{coef } ( \junul , \iunuk ) ] (R( \mu )) 
\end{equation}
for every $k,l \geq 1$ and $1 \leq \iunuk , \junul \leq n.$

$\ $

The proof of Proposition 3.8 comes out immediately by using 
the moment-cumulant formula, see e.g. \cite{S2}, Section 2.4.

$\ $

$\ $

$\ $

\setcounter{section}{4}
{\large\bf 4. A basic combinatorial remark} 

$\ $

{\bf 4.1 Definition} Let $n$ be a positive integer, let $\sigma$
\setcounter{equation}{0}
be in $NC(2n),$ and consider the permutation $Perm ( \sigma ) 
\in {\cal S}_{2n}$ associated to $\sigma$ as in Section 2.6.
We will say that $\sigma$ is {\em parity-alternating} (respectively
{\em parity-preserving}) if the difference $i-[Perm( \sigma )] (i)
\in \Z$ is odd (respectively even) for every $1 \leq i \leq 2n.$
We will use the notations
\begin{equation}
\left\{  \begin{array}{l}
{ \{ \sigma \in NC(2n) \ | \ \sigma \mbox{ parity-alternating } \}
\ \egdef \ NC_{p-alt} (2n), }  \\  
{ \{ \sigma \in NC(2n) \ | \ \sigma \mbox{ parity-preserving} \}
\ \egdef \ NC_{p-prsv} (2n). }    
\end{array}  \right.
\end{equation}

$\ $

{\bf 4.2 Remark} It is clear that a partition $\sigma \in NC(2n)$
is parity-preserving if and only if each of its blocks either is
contained in $\{ 1,3, \ldots , 2n-1 \}$ or is contained in
$\{ 2,4, \ldots , 2n \} .$ For parity-alternating partitions,
an equivalent description is:
\begin{equation}
\sigma \in NC_{p-alt} (2n) \ \Leftrightarrow \ \{ \mbox{ every block
of $\sigma$ has an even number of elements} \} .
\end{equation}
Implication ``$\Rightarrow$'' in (4.2) is immediate, while
``$\Leftarrow$'' uses the fact that if $B$ is a block of $\sigma$
and if $i<j$ in $B$ are such that 
$\{ i+1,  \ldots ,j-1 \} \cap B = \emptyset ,$ then the interval
$\{ i+1,  \ldots ,j-1 \}$ is a union of other blocks of $\sigma$
(hence the named interval has an even number of elements,
and hence $i$ and $j$ have different parities).

$\ $

{\bf 4.3 Remark} The Kreweras complementation map $K:NC(2n)
\rightarrow NC(2n)$ puts in bijection parity-alternating and 
parity-preserving partitions; we have in fact both equalities:
\begin{equation}
\left\{  \begin{array}{l}
{ K( NC_{p-alt} (2n) ) \ = \ NC_{p-prsv} (2n),  }   \\
{ K( NC_{p-prsv} (2n) ) \ = \ NC_{p-alt} (2n).  }   
\end{array}  \right.
\end{equation}
Indeed, denoting the cycle $\ccycle$ by $\gamma$, we know
(Eqn.(2.4)) that $Perm(K( \sigma )) = Perm( \sigma )^{-1} \gamma ,$
$\sigma \in NC(2n);$ this equality and the fact that
$\gamma$ itself is parity-alternating
imply together the inclusion ``$\subseteq$'' of both Equations (4.3).
But then, since $K$ is one-to-one, we get that both inequalities
between $| NC_{p-alt} (2n) |$ and $| NC_{p-prsv} (2n) |$ are holding;
hence $| NC_{p-alt} (2n) | = | NC_{p-prsv} (2n) |,$ and in (4.3)
we must really have equalities.

$\ $

The combinatorial remark mentioned in the title of the section is
the following.

$\ $

{\bf 4.4 Proposition} Let $n$ be a positive integer. In order to
distinguish the notations, we will write $\kd$ for the Kreweras
complementation map on $NC(2n)$ (while the Kreweras map on
$NC(n)$ will be simply denoted by $K$).

$1^{o}$ We have a canonical bijection 
\begin{equation}
\{ ( \pi , \rho ) \ | \ \pi , \rho \in NC(n), \ \pi \leq \rho \ \}
\ni ( \pi , \rho ) \ \longrightarrow  \ \sigma \in NC_{p-prsv} (2n),
\end{equation}
where $\sigma$ of (4.4) is determined as follows: the restriction
$\sigma | \{ 2,4, \ldots , 2n \}$ is obtained by transporting $\pi$
from $\nn$ to $\{ 2,4, \ldots , 2n \}$, while the restriction
$\sigma | \{ 1,3, \ldots , 2n-1 \}$ is obtained by transporting 
$K^{-1} ( \rho )$ from $\nn$ to $\{ 1,3, \ldots , 2n-1 \}$. (In
terms of equivalence relations associated to partitions, we have
$2i \ecsigma 2j \Leftrightarrow i \ecpi j$ and
$2i-1 \ecsigma 2j-1 \Leftrightarrow i \eckinvrho j$, for every
$1 \leq i,j \leq n.)$

$2^{o}$ If $\pi \leq \rho$ in $NC(n)$ and $\sigma \in NC_{p-prsv} (2n)$
are as above, then the 
relative Kreweras complement $K_{\rho} ( \pi ) \in NC(n)$
(described in Section 2.5) is given by the formula:
\begin{equation}
K_{\rho} ( \pi ) \ = \ \left\{ \ \{ \frac{b}{2} \ | \ b \in 
B \cap \{ 2,4, \ldots ,2n \} \} \ ; \ B
\mbox{ block of } \kd ( \sigma ) \right\}  .
\end{equation}

$\ $

{\bf Proof} $1^{o}$ is an immediate consequence of Proposition 2.4.

For $2^{o}$, we will compare the permutations associated to the 
two partitions appearing in (4.5) (the right-hand side of (4.5), 
which could be written as $\frac{1}{2} [ \kd ( \sigma ) |
\{ 2,4, \ldots , 2n \} ]$, is also in $NC(n);$ this is because the
naturally defined operation of restriction preserves the quality
of a partition of being non-crossing - easy verification).

Let us denote the permutations $Perm( \pi )$ and $Perm( \rho )$
by $\alpha$ and $\beta$, respectively. Then we know (Eqn.(2.5))
that $Perm( K_{\rho} ( \pi )) = \alpha^{-1} \beta .$ On the other
hand, let us denote the partition in the right-hand side of (4.5)
by $\theta$. We have that $\kd ( \sigma ) \in NC_{p-alt} (2n)$
(because $\sigma \in NC_{p-prsv} (2n)$ and by Remark 4.3), and
this immediately implies the following formula, expressing
$Perm ( \theta )$ in terms of $Perm( \kd ( \sigma )):$
\begin{equation}
2 [ Perm( \theta ) ] (i) \ = \ \left( Perm( \kd ( \sigma ))
\right)  ^{2} (2i), \ \ 1 \leq i \leq n.
\end{equation}
Hence what we have to do is calculate $Perm( \kd ( \sigma ))$
in terms of $\alpha$ and $\beta$, and verify that the right-hand
side of (4.6) equals $2 \alpha^{-1} ( \beta (i)).$

Now, the definition of $\sigma$ in terms of $\pi$ and $\rho$ is
converted at the level of the associated permutations by saying
that $Perm ( \sigma ) \in {\cal S}_{2n}$ acts by:
\[
\left\{  \begin{array}{l}
{ [Perm( \sigma )] (2i) \ = \ 2[ Perm( \pi )] (i) \ = \ 
2 \alpha (i), } \\
{ [Perm( \sigma )] (2i-1) \ = \ 2[ Perm( K^{-1} ( \rho ) )] (i) -1, }
\end{array}  \ \ 1 \leq i \leq n.  \right.
\]
It is more convenient to look at $Perm( \sigma )^{-1} ,$ which is
thus given by 
\begin{equation}
\left\{  \begin{array}{l}
{ [Perm( \sigma )]^{-1} (2i) \ = \ 2 \alpha^{-1} (i), } \\
{ [Perm( \sigma )]^{-1} (2i-1) \ = \ 2[ Perm( K^{-1} ( \rho ) )]^{-1}
(i) -1, }
\end{array}  \ \ 1 \leq i \leq n.  \right.
\end{equation}

Denoting the cycles $\cycle \in {\cal S}_{n}$ and 
$\ccycle \in {\cal S}_{2n}$ by $\gamma_{n}$ and $\gamma_{2n}$,
respectively, we get (by appropriately substituting in (2.4))
that $Perm( K^{-1} ( \rho ))^{-1} = \beta \gamma_{n}^{-1}$
and that $Perm( \kd ( \sigma )) = Perm( \sigma )^{-1} \gamma_{2n} .$
Hence for every $1 \leq i \leq n:$
\begin{equation}
[Perm( \kd ( \sigma )) ] (2i-1) \ = \ 
[ Perm( \sigma )^{-1} \gamma_{2n} ] (2i-1) \ = \
Perm( \sigma )^{-1} (2i) \ \stackrel{(4.7)}{=} \ 2 \alpha^{-1} (i);
\end{equation}
and for every $1 \leq i \leq n-1$, a similar calculation leads to
\begin{equation}
[Perm( \kd ( \sigma )) ] (2i) \ = \ 2 \beta (i) -1.
\end{equation}
Equation (4.9) must also hold for $i=n,$ because
$Perm( \kd ( \sigma ))$ is a bijection. From (4.8) and (4.9) it is
easily verified that the right-hand side of (4.6) is indeed
$2 \alpha^{-1} ( \beta (i)),$ and this concludes the proof.
{\bf QED}  

$\ $

{\bf 4.5 Corollary} There exists a bijection between
$\{ ( \pi , \rho ) \ | \ \pi , \rho \in NC(n),$
$\pi \leq \rho \}$ and $NC_{p-alt} (2n)$, such that, for
$( \pi , \rho )$ in the first set corresponding to $\tau$
in the second:

(i) $\pi$ and $\tau$ have the same number, say $k$, of blocks,
\newline
and moreover,

(ii) we can write $\pi = \{ A_{1}, \ldots , A_{k} \}$ and
$\tau  = \{ B_{1}, \ldots , B_{k} \}$ in such a way that
$|A_{j}| = \frac{1}{2} |B_{j}|$ for every $1 \leq j \leq k.$

$\ $

(Note: The less involved fact that $\{ ( \pi , \rho ) \ |$
$\pi, \rho \in NC(n), \ \pi \leq \rho \}$ and $NC_{p-alt} (2n)$
have the same number  - namely $(3n)!/n!(2n+1)!$ - of elements,
was well-known, see e.g. \cite{E} or Section 2 of \cite{S1}.)

$\ $

{\bf Proof} The desired bijection is the diagonal arrow of the
diagram
\begin{equation}
\begin{array}{cccccccc} 
{ \{ ( \pi , \rho ) \ | \ \pi , \rho \in NC(n), \ \pi \leq \rho \} }  &
{ -- }  &  { - \rightarrow  }  &  {  NC_{p-prsv} (2n) } &
{ - }  & \stackrel{\kd}{--} &  { \rightarrow  }  &  { NC_{p-alt} (2n) } \\
{ | }        &  &  &  &  &  &  &  \\
{ | }        &  &  &  &  &  &  &  \\
{ | }        &  &  &  &  &  &  &  \\
{ | }        &  &  &  &  &  &  &  \\
{ \downarrow }        &  &  &  &  &  &  &  \\
{ \{ ( \pi , \rho ) \ | \ \pi , \rho \in NC(n), \ \pi \leq \rho \} }  &
  &  &  &  &  &  &
\end{array}
\end{equation}
where the leftmost horizontal arrow is the bijection described
in Proposition 4.4.$1^{o}$, and the vertical arrow is the bijection
$( \pi , \rho ) \rightarrow ( K_{\rho} ( \pi ), \rho ).$ Indeed,
if $( \pi , \rho ) \rightarrow \tau$ by the diagonal arrow of (4.10),
then from Proposition 4.4.$2^{o}$ it follows that $\pi$ =  
$\left\{ \ \{ \frac{b}{2} \ | \ b \in B \cap \{ 2,4, \ldots ,2n 
\} \} \ ; \right.$
\newline
$\left. B \mbox{ block of } \tau \right\} ;$ but since $\tau$ is 
parity-alternating, it is clear that 
$| \{ B \cap \{ 2,4, \ldots ,2n \} | = |B|/2$ for every block
$B$ of $\tau .$ {\bf QED}

$\ $

We now pass to the proof of the Theorem 1.14. We use the same line
as for the mentioned Application 1.10 of \cite{NS2}, which relies 
on the following

$\ $

{\bf 4.6 Freeness Criterion:} Let $\ncps$ be a $\noncom$, such that 
$\varphi$ is a trace, let $\CC \subseteq \A$ be a unital subalgebra, 
and let $\X \subseteq \A$ be a non-void subset. Assume that for every 
$m \geq 1$ and $c_1 , \ldots , c_m \in {\cal C}$, 
$x_1 , \ldots , x_m \in \X$, it is true that
\begin{equation}
\varphi ( c_1 x_1 c_2 x_2 \cdots c_m x_m ) \ =
\ [ \mbox{coef } (1, \ldots , m) ] 
(R( \mu_{ c_1 , \ldots , c_m } ) \freestar 
M( \mu_{ x_1 , \ldots , x_m }))
\end{equation}
(where, according to the notations set in Section 3.2, the right-hand 
side of (4.11) denotes the coefficient of $z_1 \cdots z_m$ in the 
formal power series $R( \mu_{ c_1 , \ldots , c_m } ) \freestar 
M( \mu_{ x_1 , \ldots , x_m })).$ Then $\X$ is free from ${\cal C}$ 
in $\ncps$.

$\ $

For the proof of 4.6, see \cite{NS2}, Proposition 4.7. In addition
to this criterion, we will use the following simple lemma.

$\ $

{\bf 4.7 Lemma} Given $m \geq 1$ and $\pi \in NC(m),$ the 
following are equivalent:

$1^{o}$ there exists a block $B$ of $\pi$ such that $|B|$ is 
odd;

$2^{o}$ there exists a block $B = \{ i_{1} < i_{2} < \cdots
< i_{k} \}$ of $\pi$ such that either $k=1$, or $k \geq 3$
is odd and all the differences $i_{2} - i_{1}, i_{3} - i_{2},
\ldots , i_{k} - i_{k-1}$ are also odd.

$\ $

{\bf Proof} We only need to show that $1^{o} \Rightarrow 2^{o}.$
We proceed by induction on $m.$ The case $m=1$ is clear, we prove
the induction step $\{ 1, \ldots , m-1 \} \Rightarrow m.$ We fix
$\pi \in NC(m)$ which satisfies $1^{o}$, and we also fix a block
$B$ of $\pi$ such that $|B|$ is odd. Clearly, we can assume that
$|B| \neq 1$ and that there exist two elements $i<j$ in $B$ such
that $\{ i+1, \ldots , j-1 \} \cap B = \emptyset$ and 
such that $j-i$ is even (otherwise we are
done). We also fix $i$ and $j.$ The interval
$\{ i+1, \ldots, j-1 \}$ is a union of blocks of $B;$ let
$\pi_{o} \in NC(j-i-1)$ be the partition obtained from
$\pi | \{ i+1, \ldots , j-1 \}$ by translation with $i$
downwards. Then $\pi_{o}$ satisfies $1^{o}$ (because $j-i-1$
is odd), hence also $2^{o}$ (by the induction hypothesis).
It only remains to pick a block of $\pi_{o}$ which satisfies 
$2^{o},$ and shift it back (with $i$ upwards) into a block
of $\pi .$  {\bf QED}

$\ $

{\bf Proof of Theorem 1.14} Let $\ncps$ and $\aunun ,b' ,b'' \in \A$
be as in the statement of the theorem. Let ${\cal C}$ denote the 
unital subalgebra of $\A$ generated by $\aunun$, and put 
$\X = \{ b'a_{1}b'', \ldots , b'a_{n}b'' \}$. The sufficient freeness
condition of 4.6 becomes here:
\begin{equation}
\varphi( c_{1} (b'a_{i_{1}}b'') c_{2} (b'a_{i_{2}}b'') \cdots 
c_{m} (b'a_{i_{m}}b'')) 
\end{equation}
\[
\ = \ [ \mbox{coef } (1, \ldots ,m) ]
\left( R( \mu_{ c_{1} , \ldots , c_{m} } ) \freestar
M( \mu_{ b'a_{i_{1}}b'' , \ldots , b'a_{i_{m}}b'' } )  \right) ,
\]
to be verified for every 
$m \geq 1, \ c_{1} , \ldots , c_{m} \in {\cal C}$
and $1 \leq i_{1} , \ldots , i_{m} \leq n.$

We start with the left-hand side of (4.12), and rewrite it as
\[
\varphi( (c_{1} b')(a_{i_{1}}b'')(c_{2} b')(a_{i_{2}}b'') \cdots 
(c_{m} b')(a_{i_{m}}b'')) 
\]
\[
= \ [ \mbox{coef } (1,2, \ldots ,2m) ]
( M( \mu_{ c_{1} b', a_{i_{1}}b'', \ldots , c_{m} b', a_{i_{m}}b'' } )
\]
\[
= \ [ \mbox{coef } (1,2, \ldots ,2m) ]
( R( \mu_{ c_{1}, a_{i_{1}}, \ldots , c_{m}, a_{i_{m}} } )
\freestar M( \mu_{ \underbrace{b',b'', \ldots , b',b''}_{2m} } ) ),
\]
(because $\{ b' ,b'' \}$ is free from $\{ c_{1} , \ldots , c_{m}, 
a_{i_{1}}, \ldots , a_{i_{m}} \}$ and by the formula (3.15))
\begin{equation}
= \ \sum_{\sigma \in NC(2m)}
[ \mbox{coef } (1,2, \ldots ,2m); \sigma ]
( R( \mu_{ c_{1}, a_{i_{1}}, \ldots , c_{m}, a_{i_{m}} } ))
\cdot
\end{equation}
\[
[ \mbox{coef } (1,2, \ldots ,2m); \kdm ( \sigma ) ]
( M( \mu_{b',b'', \ldots , b',b''} ) )
\]
(by the combinatorial definition of $\freestar$ in 3.3).

Now, let us remark that the quantity 
$[ \mbox{coef } (1,2, \ldots ,2m); \kdm ( \sigma ) ]
( M( \mu_{b',b'', \ldots , b',b''} ) )$ appearing in (4.13) vanishes
whenever $\kdm ( \sigma )$ is not parity-alternating. Indeed, if
$\kdm ( \sigma ) \not\in NC_{p-alt} (2m),$ then $\kdm ( \sigma )$
must have at least one block $B$ with $|B|$ odd (by Remark 4.2),
hence also at least one block $B$ with the property stated in
$4.7.2^{o}$ (by the Lemma 4.7). But for any block $B$ with the latter
property we have $[ \mbox{coef } (1,2, \ldots ,2m)| B ]
( M( \mu_{b',b'', \ldots , b',b''} ) )$ =0, because $(b',b'')$ is
diagonally balanced; this implies that
$[ \mbox{coef } (1,2, \ldots ,2m); \kdm ( \sigma ) ]$
\newline
$( M( \mu_{b',b'', \ldots , b',b''} ) )$ =0, by Eqn.(3.4).

By also invoking the Remark 4.3, we hence see that in (4.13) we can 
in fact sum only over $\sigma \in NC_{p-prsv} (2m)$ (the contribution 
of a $\sigma \not\in NC_{p-prsv} (2m)$ is always zero).

But then we can perform in (4.13) the ``substitution'' indicated
by the bijection of Proposition 4.4.$1^{o}$, i.e. we can pass to 
a sum indexed by $\{ ( \pi , \rho ) \ | \ \pi , \rho \in NC(m),$
$\pi \leq \rho \} .$ Of course, in order to do this we need to
rewrite the general summand of (4.13) not in terms of 
$\sigma \in NC_{p-prsv} (2m),$ but in terms of the pair 
$( \pi , \rho )$ corresponding to it. We leave it to the reader
to verify that 

(a) the part $[ \mbox{coef } (1,2, \ldots ,2m); \sigma ]
( R( \mu_{ c_{1}, a_{i_{1}}, \ldots , c_{m}, a_{i_{m}} } ))$ of
the general summand gets converted into the product
\[
[ \mbox{coef } (1, \ldots ,m); \pi ]
( R( \mu_{ a_{i_{1}}, \ldots , a_{i_{m}} } )) \cdot
[ \mbox{coef } (1, \ldots ,m); K^{-1} ( \rho ) ]
( R( \mu_{ c_{1}, \ldots , c_{m} } ));
\]
and

(b) due to the description of $K_{\rho} ( \pi )$ in
Proposition 4.4.$2^{o}$, the part 
$[ \mbox{coef } (1,2, \ldots ,2m);$
\newline
$\kdm ( \sigma ) ]$
$( M( \mu_{ \underbrace{b',b'', \ldots , b',b''}_{2m} } ) )$ 
of the general summand gets converted into the coefficient
\newline
$[ \mbox{coef } (1, \ldots ,m); K_{ \rho} ( \pi ) ]
( M( \mu_{ \underbrace{b''b', \ldots , b''b'}_{m} } ) ).$
(Indeed, if we write $\kdm ( \sigma )$ =
$\{ B_{1} , \ldots , B_{k} \}$ and $K_{\rho} ( \pi )$ =
$\{ A_{1} , \ldots , A_{k} \}$ in such a way that $|B_{j}| =
2|A_{j}|,$ $1 \leq j \leq k,$ then the both $[ \mbox{coef } \ldots ]$
quantities claimed to be equal are seen to be just
$\prod_{j=1}^{k} \varphi ( (b''b')^{|A_{j}|} ).$ We are of course
using here the fact that $\varphi ( (b'b'')^{l} )$ =
$\varphi ( (b''b')^{l} )$ for every $l \geq 0,$ which holds because 
$\varphi$ is a trace.)

The conclusion of (a) and (b) above is that (4.13) can be
continued with:
\begin{equation}
\sum_{ \begin{array}{c}
{\scriptstyle \pi , \rho \in NC(n)} \\
{\scriptstyle such \ that}  \\
{\scriptstyle \pi \leq \rho }
\end{array}  }
[ \mbox{coef } (1, \ldots ,m); K^{-1} ( \rho ) ]
( R( \mu_{ c_{1}, \ldots , c_{m} } )) \cdot
[ \mbox{coef } (1, \ldots ,m); \pi ]
( R( \mu_{ a_{i_{1}}, \ldots , a_{i_{m}} } ))  \cdot
\end{equation}
\[
[ \mbox{coef } (1, \ldots ,m); K_{ \rho} ( \pi ) ]
( M( \mu_{ b''b', \ldots , b''b' } ) ) .
\]
It is not hard to check that the quantity in (4.14) is exactly
\begin{equation}
[ \mbox{coef } (1, \ldots ,m) ]
( R( \mu_{ c_{1}, \ldots , c_{m} } ) \freestar
R( \mu_{ a_{i_{1}}, \ldots , a_{i_{m}} } )  \freestar
M( \mu_{ b''b' , \ldots , b''b' } ) ) 
\end{equation}
(indeed, one has just to look at the combinatorial formula for the
$\freestar$-product of three series $f,g,h,$ calculated in the 
order $f \ \freestar \ ( g \ \freestar \ h))$ - compare for instance 
to \cite{NS2}, Eqn.(3.7)).

Finally, since $b''b'$ is free from $\{ a_{i_{1}}, \ldots , a_{i_{m}}
\}$ we have (by applying again (3.15)) that
$R( \mu_{ a_{i_{1}}, \ldots , a_{i_{m}} } )  \freestar
M( \mu_{ b''b' , \ldots , b''b' } )$ =
$M( \mu_{ a_{i_{1}} b''b' , \ldots , a_{i_{m}} b''b' } )$; the latter
series can also be written as
$M( \mu_{ b' a_{i_{1}} b'', \ldots , b' a_{i_{m}} b'' } )$, because 
we are working in a tracial non-commutative probability space. We thus 
arrive to the fact that (4.15) coincides with the right-hand side of 
(4.12), and this concludes the proof. {\bf QED} 

$\ $

$\ $

$\ $

\setcounter{section}{5}
{\large\bf 5. Diagonally balanced pairs} 

$\ $

We first show that Definition 1.12 could have been equally well
\setcounter{equation}{0}
formulated in terms of the free cumulants of $a_{1}$ and
$a_{2}$ (rather than their moments). 

$\ $

{\bf 5.1 Proposition} Let $\ncps$ be a $\noncom$, and let 
$\aunudoi$ be in $\A .$ The pair $( \aunudoi )$ is diagonally
balanced if and only if the coefficients of
$\underbrace{z_{1} z_{2} \cdots z_{1} z_{2} z_{1} }_{2n+1}$ and
$\underbrace{z_{2} z_{1} \cdots z_{2} z_{1} z_{2} }_{2n+1}$ in
$[ R( \mu_{\aunudoi} ) ] ( \zunudoi )$ are equal to zero for
every $n \geq 0.$

$\ $

{\bf Proof} This is an immediate consequence of the Lemma 4.7
and the Equations (3.10), (3.11) connecting the $M$- and $R$-series of
$\mu_{\aunudoi} .$ For instance, when proving the ``$\Rightarrow$''
part, we write
\begin{equation}
[ \mbox{coef } (1,2, \ldots ,1,2,1) ] (R( \mu_{\aunudoi} )) \ = \ 
[ \mbox{coef } (1,2, \ldots ,1,2,1) ] (M( \mu_{\aunudoi} ) 
\freestar Moeb )
\end{equation}
\[
= \ \sum_{\pi \in NC(2n+1)}
[ \mbox{coef } (1,2, \ldots ,1,2,1) ; \pi ] (M( \mu_{\aunudoi} )) 
\cdot [ \mbox{coef } (1,2, \ldots ,1,2,1) ; K( \pi ) ] ( Moeb ) 
\]
(and a similar formula for
$[ \mbox{coef } (2,1, \ldots ,2,1,2) ] (R( \mu_{\aunudoi} ))$ ).
Then we note that in the sum (5.1) every term is actually
vanishing; indeed, the Lemma 4.7 and the hypothesis that
$( \aunudoi )$ is a diagonally balanced pair imply together that 
$[ \mbox{coef } (1,2, \ldots ,1,2,1) ; \pi ] (M( \mu_{\aunudoi} ))$
= 0 for every $\pi \in NC(2n+1).$ {\bf QED}

$\ $

{\bf 5.2 Remark} It is useful to record here that if
$( \aunudoi )$ is a diagonally balanced pair in $\ncps$, where
$\varphi$ is a trace, then the coefficient of 
$z_{i_{1}} z_{i_{2}} \cdots z_{i_{2n+1}}$ 
in $[ R( \mu_{ \aunudoi } ) ] ( \zunudoi )$ is also vanishing
whenever $(i_{1}, i_{2}, \ldots , i_{2n+1})$ is a cyclic permutation
of $(1,2, \ldots ,1,2,1)$ or of $(2,1, \ldots , 2,1,2).$ This follows
from Propositions 5.1 and 3.8.

$\ $

Since (as it is clear from 1.3 and 5.1) every $R$-diagonal pair is 
diagonally balanced, the next result contains 
the Proposition 1.7 as a particular case. 

$\ $

{\bf 5.3 Proposition} Let $\ncps$ be a $\noncom$ such that 
$\varphi$ is a trace, and let $( \aunudoi )$ be a diagonally
balanced pair in $\ncps .$ For every $k \geq 1$ we denote by
$\alpha_{k}$ the coefficient of $(z_{1} z_{2} )^{k}$
(equivalently, of $(z_{2} z_{1} )^{k}$ ) in the series
$[ R( \mu_{\aunudoi} ) ] ( \zunudoi ),$ and we consider the series
of one variable $f(z) = \sum_{k=1}^{\infty} \alpha_{k} z^{k} .$
Then $f = R( \mu_{a_{1} a_{2}}  ) \freestar Moeb,$ where $Moeb$
is the Moebius series, as in Eqn.(1.9) (and of course
$\mu_{a_{1} a_{2} } : \C [X] \rightarrow \C$ denotes the
distribution of the product $a_{1} \cdot a_{2} ).$

$\ $

{\bf Proof} For every $n \geq 1$ we have:
\[
\varphi ( (a_{1}a_{2} )^{n} ) \ = \ 
[ \mbox{coef } ( \underbrace{1,2, \ldots , 1,2}_{2n} ) ] 
(M( \mu_{\aunudoi} )) 
\]
\begin{equation}
\stackrel{(3.12)}{=}  \ \sum_{\tau \in NC(2n)}
[ \mbox{coef } ( \underbrace{1,2, \ldots , 1,2}_{2n} ) ; \tau ] 
(R( \mu_{\aunudoi} ) ). 
\end{equation}
If the partition $\tau \in NC(2n)$ appearing in (5.2) is not in
$NC_{p-alt}(2n),$ then it has at least one block $B$ with 
$|B|$ odd (by Remark 4.2), hence it also has a block $B$
satisfying condition $2^{o}$ of Lemma 4.7. This and the description
of the diagonally balanced pair $( \aunudoi )$ in terms of
cumulants (Proposition 5.1) imply together that 
$[ \mbox{coef } ( \underbrace{1,2, \ldots , 1,2}_{2n} ) ; \tau ] 
(R( \mu_{\aunudoi} ) ) =0.$ On the other hand, if $\tau =
\{ B_{1} , \ldots , B_{h} \}$ is in $NC_{p-alt} (2n),$ then by
the very definition of $NC_{p-alt} (2n)$ we have that 
\[
[ \mbox{coef } ( \underbrace{1,2, \ldots , 1,2}_{2n} ) ; \tau ] 
(R( \mu_{\aunudoi} ) ) \ = \ \alpha_{|B_{1}|/2} \cdots
\alpha_{|B_{h}|/2}
\]
(where the $\alpha$'s are as in the statement of the proposition).
This argument shows that the sum in (5.2) is equal to:
\begin{equation}
\sum_{ \begin{array}{c}
{\scriptstyle \tau \in NC_{p-alt} (2n) } \\
{\scriptstyle \tau = \{ B_{1} , \ldots , B_{h} \}  }
\end{array}  }
 \ \alpha_{|B_{1}|/2} \cdots \alpha_{|B_{h}|/2} .
\end{equation}
We next transform the sum in (5.3) by using the bijection between
$NC_{p-alt} (2n)$ and $\{ ( \pi , \rho ) \ |$
$\pi , \rho \in NC(n),$ $\pi \leq \rho \}$ indicated in Corollary 4.5;
we obtain:
\[
\sum_{ \begin{array}{c}
{\scriptstyle \pi , \rho \in NC(n) } \\
{\scriptstyle such \ that \ \pi \leq \rho ,} \\
{\scriptstyle \pi = \{ A_{1} , \ldots , A_{h} \}  }
\end{array}  }
 \ \alpha_{|A_{1}|} \cdots \alpha_{|A_{h}|} 
\]
\begin{equation}
= \ \sum_{ \begin{array}{c}
{\scriptstyle \pi \in NC(n) } \\
{\scriptstyle \pi = \{ A_{1} , \ldots , A_{h} \}  }
\end{array}  }
 \ \alpha_{|A_{1}|} \cdots \alpha_{|A_{h}|} \cdot 
\mbox{ card } \{ \rho \in NC(n) \ | \ \rho \geq \pi \} .
\end{equation}
Now, for a given $\pi = \{ A_{1}, \ldots , A_{h} \} \in NC(n),$
the product $\alpha_{|A_{1}|} \cdots \alpha_{|A_{h}|}$
is just $[ \mbox{coef } (n); \pi ] (f)$ (in the sense of the 
notations in 3.2.$4^{o}$), while on the other hand it is an easy
exercise to identify
$\mbox{ card } \{ \rho \in NC(n) \ | \ \rho \geq \pi \}$ as
$[ \mbox{coef } (n) ; K( \pi )] (Zeta \freestar Zeta).$ Hence
(5.4) can be continued with:
\[
\sum_{\pi \in NC(n)} 
[ \mbox{coef } (n); \pi ] (f) \cdot
[ \mbox{coef } (n) ; K( \pi )] (Zeta \freestar Zeta) \ = \ 
[ \mbox{coef } (n) ] (f \freestar Zeta \freestar Zeta).
\]

But the expression we had started with was 
$\varphi ( (a_{1} a_{2} )^{n} )$, i.e. 
$[ \mbox{coef } (n) ] (M ( \mu_{ a_{1} \cdot a_{2} } )).$ We have in other
words proved the equality
\begin{equation}
M( \mu_{a_{1} \cdot a_{2}} ) \ = \ f \freestar Zeta \freestar Zeta
\end{equation}
(since the two series have identical coefficients). It only remains
to take the $\freestar$ product with $Moeb \freestar Moeb$ on both
sides of (5.5), and take into account that $Moeb$ is the inverse
of $Zeta$ under $\freestar$, and that $M( \mu_{a_{1} \cdot a_{2}} )
\freestar Moeb = R( \mu_{a_{1} \cdot a_{2}} ).$ {\bf QED}

$\ $

$\ $

$\ $

\setcounter{section}{6}
{\large\bf 6. The line of proof for the Theorem 1.5} 

$\ $

We begin by proving the particular case when the pair 
\setcounter{equation}{0}
$\{ p_{1} , p_{2} \}$ in Theorem 1.5 is of the form $\{ p,1 \} ;$ that is:

$\ $

{\bf 6.1 Proposition} Let $\ncps$ be a $\noncom$, such that
$\varphi$ is a trace, and let $a_{1} , a_{2} , p \in \A$ be such 
that $( \aunudoi )$ is an $R$-diagonal pair, and $p$ is free
from $\{ \aunudoi \} .$ Then $( a_{1} p, a_{2} )$ is also an
$R$-diagonal pair.

$\ $

{\bf Proof} We fix $m \geq 1$ and a string 
$\ee = ( \lunum ) \in { \{ 1,2 \} }^{m}$ for some $m \geq 1,$
such that $[ \mbox{coef } ( \lunum ) ] (R( \mu_{a_{1}p,a_{2}}))
\neq 0.$ Our task is to show that $m$ is even, and that the
string $\ee$ is either $( \underbrace{1,2,1,2, \ldots ,1,2}_{m} )$
or $( \underbrace{2,1,2,1, \ldots ,2,1}_{m} ).$ (This would verify
Condition (ii) in Definition 1.3; note that Condition (i) of 1.3
is automatically verified, since it is assumed that $\varphi$ is
a trace.)

The couples $( \aunudoi )$ and $(p,1)$ are free, hence we can use
Theorem 3.7 to infer that $R( \mu_{a_{1}p,a_{2}} ) = 
R( \mu_{\aunudoi} ) \freestar R( \mu_{p,1} );$ the combinatorial 
definition of $\freestar$ applied to this situation yields then
the formula:
\begin{equation}
[ \mbox{coef } ( \lunum )] (R( \mu_{a_{1}p,a_{2}} )) \ = \ 
\sum_{\tau \in NC(m)} [ \mbox{coef } ( \lunum ) ; \tau ] 
(R( \mu_{a_{1},a_{2}} ))  \cdot
\end{equation}
\[
[ \mbox{coef } ( \lunum ) ; K( \tau )] (R( \mu_{p,1} )). 
\]
Since the left-hand side of (6.1) is assumed to be non-zero,
we can choose (and fix) a partition $\tau \in NC(m)$ such that
\begin{equation}
\left\{  \begin{array}{l}
{ [ \mbox{coef } ( \lunum ) ; \tau ] 
(R( \mu_{a_{1},a_{2}} )) \neq 0 } \\ 
{ [ \mbox{coef } ( \lunum ) ; K( \tau )] 
(R( \mu_{p,1} )) \neq 0. }
\end{array}  \right.
\end{equation}

The first condition (6.2) together with the hypothesis that
$( \aunudoi )$ is an $R$-diagonal pair imply together that 
for every block $B= \{ i_{1} < i_{2} < \cdots < i_{k} \}$
of $\tau :$ $k$ is even, and $l_{i_{1}} \neq l_{i_{2}},$
$l_{i_{2}} \neq l_{i_{3}}, \ldots , l_{i_{k-1}} \neq l_{i_{k}}.$
This already implies that $m$ is even, and that among $\lunum$
there are $m/2$ occurrences of 1 and $m/2$ occurrences of 2.

In order to get the interpretation for the second condition
(6.2), we note that, since $p$ is always free from 1, we have
(by Theorem 3.6, in fact its particular case stated in Eqn.(1.4)):
\begin{equation}
[R( \mu_{p,1} )] ( \zunudoi ) \ = \ 
[R( \mu_{p} )] ( z_{1} ) + [R( \mu_{1} )] ( z_{2} ) \ = \ 
[R( \mu_{p} )] ( z_{1} ) + z_{2} .
\end{equation}
By taking this into account, we see that if the second
condition (6.2) is satisfied, then $\{ i \}$ has to be a
one-element block of $K( \tau )$ whenever we have $l_{i} =2.$
A quick look at how the Kreweras complement 
$K:NC(m) \rightarrow NC(m)$ is defined shows that $\{ i \}$
is a one-element block of $K( \tau )$ if and only if $i$ and
$i+1$ are in the same block of $\tau$ (where if $i=m$ we read
``1'' instead of ``$i+1$''). Hence whenever we have $l_{i} =2,$
we also have that $i$ and $i+1$ lie in the same block of $\tau$,
and then from the discussion in the preceding paragraph we deduce
that $l_{i+1} =1.$

We have thus proved that: $m$ is even; among $\lunum$ there are
$m/2$ of 1 and $m/2$ of 2; and $l_{i}=2 \Rightarrow l_{i+1} =1,$
for $1 \leq i \leq m-1,$ and $l_{m} =2 \Rightarrow l_{1} =1.$
We leave it as an elementary exercise to the reader to verify that
a string $( \lunum )$ with all these properties can only be
$(1,2,1,2, \ldots ,1,2)$ or $(2,1,2,1, \ldots ,2,1).$ {\bf QED}

$\ $

{\bf 6.2 Remark} With exactly the same argument, we could have
shown that $( \aunudoi )$ $R$-diagonal $\Rightarrow$
$( a_{1}, pa_{2} )$ $R$-diagonal (i.e. the particular case
$\{ \punudoi \} = \{ 1,p \}$ of Theorem 1.5). Since (obviously)
$R$-diagonality is not affected by reversing the order of the
elements of a pair, it follows that in Proposition 6.1 the
element $p$ could have been in fact inserted anywhere (left or
right of either $a_{1}$ or $a_{2}).$ We take the occasion to 
mention here that (as the reader can easily verify by 
experimenting on coefficients of small length) the Theorem 1.5
itself becomes false if in its statement the pair
$( a_{1}p_{1} , p_{2}a_{2} )$ is replaced for instance with
$( a_{1}p_{1} , a_{2}p_{2} ).$

$\ $

In view of what has been proved up to now, the Theorem 1.5 is
equivalent to the following statement.

$\ $

{\bf 6.3 Proposition} Let $\ncps$ be a $\noncom$, such that 
$\varphi$ is a trace, and let $\aunudoi , \punudoi \in \A$
be such that $( \aunudoi )$ is an $R$-diagonal pair, and
such that $\{ \punudoi \}$ is free from $\{ \aunudoi \} .$
Then $\mu_{a_{1}p_{1} , p_{2}a_{2}} = 
\mu_{a_{1}p_{1}p_{2}, a_{2}} .$

$\ $

Indeed, Theorem 1.5 follows from Propositions 6.1 and 6.3,
while conversely, Proposition 6.3 is implied by Theorem 1.5
and Corollary 1.8.

Proposition 6.3 simply says that if we put (in the notations
used there) $x_{1} = a_{1}p_{1}, x_{2} = p_{2}a_{2},$
$y_{1} = a_{1}p_{1}p_{2}, y_{2} = a_{2},$ then we have
\begin{equation}
\varphi ( x_{l_{1}} x_{l_{2}} \cdots x_{l_{m}} ) \ = \ 
\varphi ( y_{l_{1}} y_{l_{2}} \cdots y_{l_{m}} ) 
\end{equation}
for every $m \geq 1$ and $\lunum \in \{ 1,2 \} .$ The point
of view we will take concerning (6.4) is the following:
when we write back $x_{1} = a_{1}p_{1}, x_{2} = p_{2}a_{2},$
$y_{1} = a_{1}p_{1}p_{2}, y_{2} = a_{2},$ then both
$x_{l_{1}} x_{l_{2}} \cdots x_{l_{m}}$ and
$y_{l_{1}} y_{l_{2}} \cdots y_{l_{m}}$ are converted into
monomials of length $2m$ in $\aunudoi, \punudoi ;$ and
moreover, both these monomials of length $2m$ can be regarded
as the result obtained by starting with
$a_{l_{1}} a_{l_{2}} \cdots a_{l_{m}}$ and inserting 
$p_{1}$'s and $p_{2}$'s in between the $a$'s, according to a 
certain pattern:

- in order to obtain $x_{l_{1}} x_{l_{2}} \cdots x_{l_{m}}$,
we insert a $p_{1}$ immediately to the right of each occurrence 
of $a_{1} ,$ and a $p_{2}$ immediately to the left of each
occurrence of $a_{2};$

- in order to obtain $y_{l_{1}} y_{l_{2}} \cdots y_{l_{m}}$,
we insert a $p_{1}p_{2}$ immediately to the right of each 
occurrence of $a_{1}.$

So in some sense what we have to do is compare the two
insertions patterns described above. Due to its length,
this will be done separately in the next section.

We only make now one simple remark: since $\varphi$ is a trace,
both monomials appearing in (6.4) can be permuted cyclically;
hence if we assume in (6.4) that $l_{1} =1,$ we are in fact
missing only the case when $l_{1} = l_{2} = \cdots = l_{m} =2.$
But the latter case can be easily settled if we note that
$\varphi ( b_{1}^{n} ) = \varphi ( b_{2}^{n} ) =0$ for every
$R$-diagonal pair $( \bunudoi )$ and for every $n \geq 1$
(this in turn is an immediate consequence of the moment-cumulant
formula (3.12), combined with the particular form of the series
$R( \mu_{\bunudoi} )).$ Indeed, if 
$l_{1} = l_{2} = \cdots = l_{m} =2,$ then the right-hand side
of (6.4) becomes $\varphi ( a_{2}^{m} )=0.$ In order to verify
that the left-hand side of (6.4) is also vanishing, we can
invoke the same argument, where we use in addition the Remark 6.2
($p_{2}a_{2}$ enters the $R$-diagonal pair 
$( a_{1} , p_{2}a_{2} ),$ hence
$ \varphi ( ( p_{2}a_{2} )^{m} ) =0).$

$\ $

$\ $

\setcounter{section}{7}
{\large\bf 7. Insertion patterns associated to a string of 
1's and 2's} 

$\ $

{\bf 7.1 Notations} Throughout this section we fix a positive 
\setcounter{equation}{0}
integer $m \geq 1,$ and a string $\ee = ( \lunum ) \in \{ 1,2
\} ^{m} .$ We make the assumption that $l_{1} =1.$ The number 
of occurrences of 1 in the string will be denoted by $n$ 
($1 \leq n \leq m).$ 

Also, for the whole section we fix a circle  of radius 1 in 
the plane, and the points $\Punum , \Qunum , \Runun$ on the circle, 
positioned as follows. First we draw $P_{1}, \ldots ,$
\newline
$P_{m}$ around the circle, 
equidistant and in clockwise order. Then we draw $\Qunum ,$
according to the following rule: if $l_{i} =1,$ we put $Q_{i}$
on the arc of circle going from $P_{i}$ to $P_{i+1}$, such that
the length of the arc from $P_{i}$ to $Q_{i}$ is 1/3 of the length
of the arc from $P_{i}$ to $P_{i+1}$ (i.e. $2 \pi /3m);$ and if 
$l_{i} =2,$ we put $Q_{i}$ on the arc of circle from $P_{i-1}$ 
to $P_{i}$, such that the length of the arc $Q_{i} P_{i}$ is 
$2 \pi /3m$ (all arcs described clockwisely). It is convenient
to view the points $\Qunum$ as colored, namely we will say that
$Q_{i}$ is red whenever $l_{i} =1,$ and that it is blue whenever
$l_{i} =2.$ (Note that $Q_{1}$ is red, since it is assumed that
$l_{1} =1.)$ Finally, we inspect the $n$ points $Q_{i}$ that 
are red, clockwisely and starting from $Q_{1},$ and we second-name
as $\Runun$, in this order (every $Q_{i}$ which is red is hence 
at the same time an $R_{j}$ for some $j).$

For example, the next figure shows how our circular picture 
looks if $m=8$ and $\ee =(1,2,2,1,1,2,1,2).$

$\ $

\vspace{5cm}

\[
\mbox{\bf Figure 2.}
\]

{\bf 7.2 The ``complementation'' maps $C_{Q}$ and $C_{R}$}
Using the circular picture associated to the string $\ee ,$
we will define two maps $\cq : NC(m) \rightarrow NC(m),$ and
$\crr : NC(m) \rightarrow NC(n).$ It is convenient to formulate
the definition in terms of equivalence relations associated to
partitions (recall from 2.1 that if $\pi$ is a partition of 
$\kk ,$ then the equivalence relation $\ecpi$ corresponding
to $\pi$ is simply ``$i \ecpi j \Leftrightarrow i,j$ belong to 
the same block of $\pi$'', $1 \leq i,j \leq k).$ Given 
$\sigma \in NC(m),$ the equivalence relations corresponding
to the partitions $\cq ( \sigma )$ of $\mm$ and 
$\crr ( \sigma )$ of $\nn$ are defined as follows:
\begin{equation}
\left\{  \begin{array}{l}
{ \mbox{ - For $1 \leq i,j \leq m$ we have $i \eccqsigma j$
if and only if there are no $1 \leq h,k \leq m$} } \\ 
{ \mbox{   such that $h \ecsigma k$ and such that the line 
segments $Q_{i} Q_{j}$ and $P_{h} P_{k}$ have }  } \\
{ \mbox{   non-void intersection;} } \\
{       }  \\ 
{ \mbox{ - For $1 \leq i,j \leq n$ we have $i \esigma j$
if and only if there are no $1 \leq h,k \leq m$} } \\ 
{ \mbox{   such that $h \ecsigma k$ and such that the line 
segments $R_{i} R_{j}$ and $P_{h} P_{k}$ have }  } \\
{ \mbox{   non-void intersection;} } 
\end{array}  \right.
\end{equation}

We leave it to the reader to verify that $\eccqsigma$ and
$\esigma$ defined in (7.1) really are equivalence relations
on $\mm$ and $\nn$, respectively, and moreover, that
$\cq ( \sigma)$ and $\crr ( \sigma )$ are indeed non-crossing
partitions. All these verifications reduce to the same geometric
argument, that can be stated as follows: let $X,Y$ be points on 
the circle, and let $Z$ be either on the circle 
or in the open disk enclosed by it, such that $\{ X,Y,Z \}
\cap \{ \Punum \} = \emptyset ;$ if the segment
$P_{h} P_{k}$ (for some $1 \leq h < k \leq m)$ intersects 
$XY$, then it must also intersect $XZ \cup YZ.$

$\ $

The relevance of the complementation maps $\cq$ and $\crr$ for
our discussion comes from the following

$\ $

{\bf 7.3 Proposition:} Let $\ncps$ be a $\noncom$, and let
$\aunudoi , p_{1},$
\newline
$p_{2} \in \A$ be such that $\{ \aunudoi \}$
is free from $\{ \punudoi \}$ in $\ncps .$ Define 
$x_{1} = a_{1}p_{1}, x_{2} = p_{2}a_{2}, y_{1} = a_{1}p_{1}p_{2},
y_{2} =a_{2}.$ Then we have:
\[
\varphi ( x_{l_{1}} x_{l_{2}} \cdots x_{l_{m}} ) \ = \ 
\]
\begin{equation}
= \ \sum_{\sigma \in NC(m)}
[ \mbox{coef } ( \lunum ); \sigma ] ( R( \mu_{\aunudoi} )) 
\cdot [ \mbox{coef } ( \lunum ); \cq ( \sigma ) ] 
( M( \mu_{\punudoi} )) 
\end{equation}
and
\[
\varphi ( y_{l_{1}} y_{l_{2}} \cdots y_{l_{m}} ) \ = \ 
\]
\begin{equation}
= \ \sum_{\sigma \in NC(m)}
[ \mbox{coef } ( \lunum ); \sigma ] ( R( \mu_{\aunudoi} )) 
\cdot [ \mbox{coef } ( n ); \crr ( \sigma ) ] 
( M( \mu_{p_{1} \cdot p_{2}} )) 
\end{equation}
(where $\ee = ( \lunum )$ is the string fixed in the beginning
of the section; the notations for the series 
$M( \mu_{\ldots} ), R( \mu_{\ldots} ),$ and for their coefficients
are as in Section 3.2).

$\ $

{\bf Proof} We will only show the proof of (7.2), the one for
(7.3) is similar. This kind of argument has been used several
times before, in connection to the Kreweras complementation
map $K$ (instead of $\cq , \crr$) - see 3.4 in \cite{S2},
3.4 in \cite{NS1}, 3.11 in \cite{NS2}.

If in $x_{l_{1}} x_{l_{2}} \cdots x_{l_{m}}$ we write back each
$x_{1}$ as $a_{1}p_{1}$ and each $x_{2}$ as $p_{2}a_{2}$, 
we obtain a monomial of length $2m$ in $\aunudoi , \punudoi$.
By looking just at the $a$'s in this monomial, we see that
they are $a_{l_{1}}, a_{l_{2}}, \ldots , a_{l_{m}},$ exactly in
this order, and placed on a certain set of positions
$I \subseteq \mmmm ,$ with $|I| =m.$ It is also clear that on
the complementary set of positions $J = \mmmm \setminus I$ of
our monomial of length $2m$ we have 
$p_{l_{1}}, p_{l_{2}}, \ldots , p_{l_{m}},$ exactly in this
order. So we can say that 
$x_{l_{1}} x_{l_{2}} \cdots x_{l_{m}}$ 
is obtained by shuffling together
$a_{l_{1}} a_{l_{2}} \cdots a_{l_{m}}$ and 
$p_{l_{1}} p_{l_{2}} \cdots p_{l_{m}}$, where the $a$'s have 
to sit on the positions indicated by $I$,
and the $p$'s have to sit on the positions indicated by 
$J = \mmmm \setminus I.$

Now, the value of $\varphi$ on the monomial of length $2m$
discussed in the preceding paragraph can be viewed as a
coefficient of length $2m$ of the series
$M( \mu_{\aunudoi , \punudoi} ).$ We write this coefficient 
as a summation over $NC(2m),$ by using the moment-cumulant 
formula (3.12). Because 
of the assumption that $\{ \aunudoi \}$ is free from
$\{ \punudoi \}$ (which implies that we have
$[R( \mu_{\aunudoi , \punudoi} )] ( \zunudoi , z_{3} , z_{4} )$
= $[R( \mu_{\aunudoi} )] ( \zunudoi )$ +
$[R( \mu_{\punudoi} )] ( z_{3} , z_{4} )),$
the summation over $NC(2m)$ we have arrived to has a lot of
vanishing terms; a partition $\theta \in NC(2m)$ can in fact bring
a non-zero contribution to the sum if and only if each block of 
$\theta$ either is contained in $I$ or is contained in $J$ (with
$I,J$ the sets of positions of the preceding paragraph).
This brings us to the formula:
\[
\varphi ( x_{l_{1}} x_{l_{2}} \cdots x_{l_{m}} ) \ =  
\]
\begin{equation}
= \ \sum_{ \begin{array}{c}
{\scriptstyle \sigma , \tau \in NC(m) } \\
{\scriptstyle I,J - compatible}
\end{array}  }
[ \mbox{coef } ( \lunum ); \sigma ] (R( \mu_{\aunudoi} ))
\cdot
[ \mbox{coef } ( \lunum ); \tau ] (R( \mu_{\punudoi} ));
\end{equation}
in (7.4), the fact that $\sigma , \tau \in NC(m)$ are 
$I,J$-compatible has the meaning that if we transport $\sigma$
from $\mm$ onto $I$ and we transport $\tau$ from $\mm$ onto 
$J = \mmmm \setminus I$, then the partition of $\mmmm$ which 
is obtained in this way is still non-crossing.

If we now look back at how the complementation map 
$\cq : NC(m) \rightarrow NC(m)$ was defined, it is immediate
that the $I,J$-compatibility of $\sigma , \tau \in NC(m)$ 
is equivalent to the condition
$\tau \leq \cq ( \sigma ).$ This means that (7.4) can be 
continued with:
\begin{equation}
\sum_{\sigma \in NC(m)}
[ \mbox{coef } ( \lunum ); \sigma ] (R( \mu_{\aunudoi} ))
\cdot \{  \sum_{ \begin{array}{c}
{\scriptstyle \tau \in NC(m), } \\
{\scriptstyle \tau \leq \cq ( \sigma )}
\end{array}  }
[ \mbox{coef } ( \lunum ); \tau ] (R( \mu_{\punudoi} ))  \ \} .
\end{equation}
Finally, we note that
\begin{equation}
\sum_{ \begin{array}{c}
{\scriptstyle \tau \in NC(m), } \\
{\scriptstyle \tau \leq \cq ( \sigma )}
\end{array}  }
[ \mbox{coef } ( \lunum ); \tau ] (R( \mu_{\punudoi} ))
\ = \ [ \mbox{coef } ( \lunum ); \cq ( \sigma ) ] 
(M( \mu_{\punudoi} )),
\end{equation}
as it follows by a repeated utilization of the moment-cumulant
formula (3.12).
Substituting (7.6) in (7.5) brings us to the right-hand side
of (7.2), and concludes the proof. {\bf QED}

$\ $

{\bf 7.4 Corollary} Let $\ee = ( \lunum )$ be the string of
1's and 2's fixed in 7.1, and recall that $1 \leq n \leq m$ is
the number of occurrences of 1 in the string. Let on the 
other hand $\ncps$ be a $\noncom$, with $\varphi$ a trace, 
and consider $\aunudoi , \punudoi \in \A$ such that 
$( \aunudoi )$ is an $R$-diagonal pair, and such that
$\{ \punudoi \}$ is free from $\{ \aunudoi \}$. Define 
$x_{1} = a_{1}p_{1}, x_{2} = p_{2}a_{2}, y_{1} = a_{1}p_{1}p_{2},
y_{2} =a_{2}.$ If it is not true that $m$ is even and $n=m/2,$ 
then $\varphi ( x_{l_{1}} x_{l_{2}} \cdots x_{l_{m}} )$ =
$\varphi ( y_{l_{1}} y_{l_{2}} \cdots y_{l_{m}} )$ = 0. 

$\ $

{\bf Proof} If it is not true that $m$ is even and $n=m/2,$
then due to the particular form of $R( \mu_{\aunudoi} )$
we have that
$[ \mbox{coef } ( \lunum ); \sigma ] (R( \mu_{\aunudoi} ))$
= 0 for every $\sigma \in NC(m).$ But then all the terms of 
the summations on the right-hand side of (7.2), (7.3) are
equal to zero. {\bf QED}

$\ $

Corollary 7.4 shows that the equality (6.4), whose proof is 
the goal of the present section, takes place
in a trivial way unless we impose the following

$\ $

{\bf 7.5 Supplementary hypothesis:} From now on, until the
end of section, we will assume that the string 
$\ee = ( \lunum )$ fixed in 7.1 contains an equal number of 1's 
and of 2's (i.e., $m$ is even and $n=m/2).$

$\ $

In this case, the quantities
$\varphi ( x_{l_{1}} x_{l_{2}} \cdots x_{l_{m}} )$ 
and $\varphi ( y_{l_{1}} y_{l_{2}} \cdots y_{l_{m}} )$  
mentioned in 7.4 are in general non-zero, and in order to verify
their equality we will need to prove some facts concerning
``$\ee$-alternating partitions''.

$\ $

{\bf 7.6 Definition} A partition $\sigma \in NC(m)$ will be
called {\em $\ee$-alternating} (where $\ee$ = $(l_{1} , \ldots ,$ 
\newline
$l_{m} )$ is the string of 7.1) if for every block 
$B= \{ i_{1} < i_{2} < \cdots < i_{k} \}$ of $\sigma$ we have that
$l_{i_{1}} \neq l_{i_{2}} , \ldots ,
l_{i_{k-1}} \neq l_{i_{k}} , l_{i_{k}} \neq l_{i_{1}} .$ 
We denote the set of all $\ee$-alternating partitions in
$NC(m)$ by $\eealt$.

$\ $

Other two (obviously equivalent) ways of stating that $\sigma \in
NC(m)$ is $\ee$-alternating are:

- in the terminology of 2.2: we have $l_{i} \neq l_{j}$ whenever
$1 \leq i < j \leq m$ belong to the same block of $\sigma ,$ and are
consecutive in that block; or

- in the terminology of 3.2.$2^{o}$: for every block $B$ of $\sigma ,$
the $|B|$-tuple $( \lunum )|B$ is of the form $(1,2,1,2, \ldots ,1,2)$
or $(2,1,2,1, \ldots ,2,1).$

Note that from 7.6 (or any of the two equivalent reformulations
stated above) it follows that every block of $\sigma \in \eealt$ has 
an even number of elements; i.e., $\eealt$ is a subset of
$NC_{p-alt} (m)$ discussed in Section 4.

$\ $

{\bf 7.7 Proposition} The complementation map
$\cq : NC(m) \rightarrow NC(m)$ sends $\eealt$ into itself.

$\ $

{\bf Proof} Consider the picture with $2m$ points $\Punum ,
\Qunum$ sitting on a circle of radius 1, which was used to 
define $\cq$ in 7.2. We denote by $\dd$ the closed disk 
enclosed by the circle. During this proof we will be particularly
interested in a certain type of convex subsets of $\dd ,$ described
as follows. Let $i_{1}, j_{1}, \ldots , i_{k}, j_{k}$
$(1 \leq k \leq m/2)$ be a family of $2k$ distinct indices in
$\mm$, with the property that for every $1 \leq h \leq k,$ all
the points $P_{i_{1}}, P_{j_{1}}, \ldots , P_{i_{h-1}}, P_{j_{h-1}},$
$P_{i_{h+1}}, P_{j_{h+1}}, \ldots , P_{i_{k}}, P_{j_{k}}$ 
lie in the same open half-plane $\sps_{h}$ determined by the line
through $P_{i_{h}}$ and $P_{j_{h}}.$ Then the set 
$\xx \egdef \dd \cap \sps_{1} \cap \cdots \cap \sps_{k}$ will be
called the {\em trimmed disk} determined by
$i_{1}, j_{1}, \ldots , i_{k}, j_{k}$. A picture of a trimmed
disk (exemplified for $m=12)$ is shown in Figure 3:

\vspace{7cm}

\[
\mbox{\bf Figure 3:} 
\mbox{ An example of trimmed disk, when $m=12.$}
\]

The trimmed disk $\xx$ determined by
$i_{1}, j_{1}, \ldots , i_{k}, j_{k}$ is a convex set (neither open
nor closed). The points  
$P_{i_{1}}, P_{j_{1}}, \ldots , P_{i_{k}}, P_{j_{k}}$ 
will be called the {\em vertices} of $\xx$ (and can be identified
as those $P_{i} , \ 1 \leq i \leq m,$ which lie in the closure of
$\xx$ but not in $\xx$ itself). The boundary of $\xx$ consists of
$2k$ ``edges''; $k$ of these edges are the line segments
$P_{i_{h}} P_{j_{h}}, \ 1 \leq h \leq k,$ and the other $k$ edges
are arcs of the circle. Note that when we travel around the 
boundary of $\xx ,$ the rectilinear and curvilinear edges alternate.

Let us now fix for the rest of the proof a partition $\sigma \in
\eealt ,$ about which we want to show that $\cq ( \sigma )$ is 
also in $\eealt .$

For every block $B$ of $\sigma$ we denote by $\hh_{B}$ the closed 
convex hull of the points $\{ P_{i} \ | \ i \in B \} ;$
$\hh_{B}$ is thus a closed convex polygon with $|B|$ vertices,
inscribed in the circle. The polygons $\hbblock$ are disjoint,
due to the fact that $\sigma$ is non-crossing. We look at the 
connected components of the complement $\dhb$. Each of these
connected components is a trimmed disk; the proof of this fact is 
most easily done by taking the polygons $\hh_{B}$ out of $\dd$
one by one, and using an induction argument.

Given $1 \leq i,j \leq m,$ we have that $i$ and $j$ are in the 
same block of $\cq ( \sigma )$ if and only if the points
$Q_{i}$ and $Q_{j}$ lie in the same connected component of
$\dhb$; this follows immediately from the definition of
the map $\cq$ in (7.1), where at ``$\Leftarrow$'' we also use
the fact that the connected components of $\dhb$ are convex.
The blocks of $\cq ( \sigma )$ are thus found by looking at the
connected components $\xx$ of $\dhb$, with the property that
$\xx \cap \{ \Qunum \} \neq \emptyset .$ Hence in order to show that
$\cq ( \sigma )$ is $\ee$-alternating, we have to consider such 
a connected component $\xx ,$ and prove that when we travel
around the boundary of $\xx$ we meet an even number of points 
$Q_{i} ,$ which are of alternating colors. (Recall from 7.1 that
the points $\Qunum$ are colored - $Q_{i}$ is red when $l_{i}=1,$
and is blue when $l_{i}=2.)$ 

We fix for the rest of the proof a connected component $\xx$
of $\dhb$, such that $\xx \cap \{ \Qunum \} \neq \emptyset .$
As remarked earlier, $\xx$ is a trimmed disk. Note that all
the curvilinear edges of $\xx$ have length $2 \pi /m$ (because
otherwise $\xx$ would also contain some of the points 
$P_{i}, \ 1 \leq i \leq m;$ but all the points $P_{i}$ are in
$\cup_{B} \hh_{B} ).$ The points $Q_{i}$
that are to be found in $X$ are all lying on the curvilinear 
edges of the boundary of $\xx .$ On the other hand, if the arc 
from $P_{i}$ to $P_{i+1}$ is such a curvilinear edge, then 
it contains zero, one or two points from $\{ \Qunum \} ,$
according to what are $l_{i}$ and $l_{i+1}:$

- if $l_{i} = l_{i+1} =1,$ then it contains $Q_{i} ;$

- if $l_{i} = l_{i+1} =2,$ then it contains $Q_{i+1} ;$

- if $l_{i} =1, l_{i+1} =2,$ then it contains  $Q_{i}$ and $Q_{i+1} ;$

- if $l_{i} =2, l_{i+1} =1,$ then it contains no point from
$\{ \Qunum \} .$

\vspace{10pt}
The conclusion of the preceding paragraph is that if we want to find 
the points from $\{ \Qunum \}$ that lie in $\xx ,$ what we have to do
is inspect the vertices of $\xx ,$ and see for each such vertex
$P_{i}$ whether the corresponding $Q_{i}$ is on a curvilinear edge of 
$\xx$ or not.

But if we are to inspect the vertices of $\xx ,$ the way we want to 
do this is by pairing them along the rectilinear edges of $\xx .$
(Recall that the rectilinear and curvilinear edges of $\xx$ are
alternating, so that the rectilinear edges contain all the vertices
of $\xx ,$ without repetitions.) So let us go around the boundary
of $\xx ,$ clockwisely, and take one by one the rectilinear edges
which occur. Always for such an edge, call it $P_{i} P_{j}$ (where
$P_{i}$ comes first in clockwise order), we have that 
$l_{i} \neq l_{j}$ - here is the place where we are using the
hypothesis $\sigma \in \eealt .$ It is easily seen that:

- if $l_{i} =1, l_{j} =2,$ then $P_{i}$ and $P_{j}$ do not produce
points $Q_{i},Q_{j}$ in $\xx ;$

- if $l_{i} =2, l_{j} =1,$ then both $Q_{i}$ and $Q_{j}$ are in
$\xx ;$ moreover, when recording clockwisely what is
$\{ \Qunum \} \cap \xx ,$ these two points $Q_{i}, Q_{j}$ will
be consecutive, of different colors, and the blue one will come
first.

Hence we can organize the points in $\{ \Qunum \} \cap \xx$ 
in pairs, such that when we travel clockwisely around the boundary
of $\xx$ the points from each pair are consecutive, of different
colors, and with the blue one coming first. But this clearly
implies that $\{ \Qunum \} \cap \xx$ has an even number of points,
and of alternating colors - as it was to be shown. {\bf QED}

$\ $

{\bf 7.8 Corollary} If $\sigma \in \eealt$, then the partitions
$\cq ( \sigma ) \in NC(m)$ and $\crr ( \sigma ) \in NC(n)$ 
defined in 7.2 have the same number of blocks, say $k;$ moreover,
we can write them as 
$\cq ( \sigma ) = \{ B_{1} , \ldots , B_{k} \}$ and 
$\crr ( \sigma ) = \{ A_{1} , \ldots , A_{k} \}$, in such a way 
that $|B_{j}| =2|A_{j}|$ for every $1 \leq j \leq k$ (recall 
that $m=2n,$ due to the supplementary hypothesis made in 7.5).

$\ $

{\bf Proof} Consider the closed convex polygons $\hbblock$
defined as in the proof of Proposition 7.7. As pointed out in
the named proof, the blocks of $\cq ( \sigma )$ are in one-to-one 
correspondence with the connected components $\xx$ of $\dhb ,$
having the property that $\xx \cap \{ \Qunum \} \neq \emptyset .$
In exactly the same way it is seen that the blocks of $\crr ( \sigma )$
are in one-to-one correspondence with the connected components 
$\xx$ of $\dhb ,$ having the property that 
$\xx \cap \{ \Runun \} \neq \emptyset .$ So the proof will be over
if we can show that for an arbitrary connected component $\xx$
of $\dhb$ we have:
\begin{equation}
\left\{  \begin{array}{l}
{ \mbox{ (a) } 
\xx \cap \{ \Qunum \} \neq \emptyset  \Leftrightarrow
\xx \cap \{ \Runun \} \neq \emptyset ; } \\
{     }  \\
{ \mbox{ (b) if the equivalent statements of (a) are holding, } } \\
{ \mbox{     then } | \xx \cap \{ \Qunum \} |=2| \xx \cap \{ \Runun \} |. } 
\end{array}  \right.
\end{equation}
But recall now that the set $\{ \Runun \}$ is nothing else than 
the set of those points from $\{ \Qunum \}$ that are colored in red.
Hence (7.7) can be restated as:
\begin{equation}
\left\{  \begin{array}{l}
{ \mbox{ (a) there are some points $Q_{i}, \ 1 \leq i \leq m$ 
in $\xx$ if and only if} } \\
{ \mbox{     there are some red points 
$Q_{i}, \ 1 \leq i \leq m$ in $\xx$;} } \\
{   }  \\
{ \mbox{ (b) if the equivalent statements of (a) are holding, 
then the total} }  \\ 
{ \mbox{     number of points $Q_{i}$ in $\xx$ is twice the number 
of red points $Q_{i}$ in $\xx$.}  } 
\end{array}  \right.
\end{equation}
But (7.8) is a clear consequence of the fact that $\cq ( \sigma )
\in \eealt ,$ proved in 7.7. {\bf QED}

$\ $

The next proposition concludes the proof of (6.4) (and hence of
Proposition 6.3 and of Theorem 1.5).

$\ $

{\bf 7.9 Proposition} Let $\ee = ( \lunum )$ be the string of
1's and 2's fixed in 7.1; recall that according to 7.5, $m$ is
even, and the number of occurrences of both 1 and 2 in the 
string is $n=m/2.$ Let on the 
other hand $\ncps$ be a $\noncom$, with $\varphi$ a trace, 
and consider $\aunudoi , \punudoi \in \A$ such that 
$( \aunudoi )$ is an $R$-diagonal pair, and such that
$\{ \punudoi \}$ is free from $\{ \aunudoi \}$. Define 
$x_{1} = a_{1}p_{1}, x_{2} = p_{2}a_{2}, 
y_{1} = a_{1}p_{1}p_{2}, y_{2} =a_{2}.$ 
Then $\varphi ( x_{l_{1}} x_{l_{2}} \cdots x_{l_{m}} )$ =
$\varphi ( y_{l_{1}} y_{l_{2}} \cdots y_{l_{m}} ).$ 

$\ $

{\bf Proof} We consider the expressions of 
$\varphi ( x_{l_{1}} x_{l_{2}} \cdots x_{l_{m}} )$ and
$\varphi ( y_{l_{1}} y_{l_{2}} \cdots y_{l_{m}} )$ 
via summations over $NC(m),$ as shown in Eqns.(7.2),(7.3);
we will prove that for every $\sigma \in NC(m),$ the terms
indexed by $\sigma$ in the two sums of (7.2) and (7.3)
coincide. If $\sigma \not\in \eealt$ this is clear, because
$[ \mbox{coef } ( \lunum ); \sigma ] (R( \mu_{\aunudoi} ))$
= 0 (due to the particular form of $R( \mu_{\aunudoi} )$),
hence the terms indexed by $\sigma$ in (7.2) and (7.3) are
both equal to zero. For $\sigma \in \eealt$, it 
suffices to show that
\begin{equation}
[ \mbox{coef } ( \lunum ); \cq ( \sigma ) ] 
(M( \mu_{\punudoi} )) \ = \ 
[ \mbox{coef } ( m/2 ); \crr ( \sigma ) ] 
(M( \mu_{p_{1} \cdot p_{2}} )) . 
\end{equation}
According to Corollary 7.8, we can write
$\cq ( \sigma ) = \{ B_{1}, \ldots , B_{k} \} ,$
$\crr ( \sigma ) = \{ A_{1}, \ldots , A_{k} \} ,$ in such a way
that $|B_{j}| = 2|A_{j}|$ for every $1 \leq j \leq k.$ On the
other hand, due to the fact that $\cq ( \sigma )$ is also
$\ee$-alternating (by Proposition 7.7), we see that 
\begin{equation}
[ \mbox{coef } ( \lunum ); \cq ( \sigma ) ] 
(M( \mu_{\punudoi} )) \ = \ 
\varphi ( (p_{1}p_{2})^{|B_{1}|/2} ) \cdots
\varphi ( (p_{1}p_{2})^{|B_{k}|/2} )
\end{equation}
(in (7.10) the fact that $\varphi$ is a trace is also used).
But it is clear that
\begin{equation}
[ \mbox{coef } ( m/2 ); \crr ( \sigma ) ] 
(M( \mu_{p_{1} \cdot p_{2}} )) \ = \ 
\varphi ( (p_{1}p_{2})^{|A_{1}|} ) \cdots
\varphi ( (p_{1}p_{2})^{|A_{k}|/2} );
\end{equation}
the right-hand sides of (7.10) and (7.11) coincide,
which establishes (7.9). {\bf QED}

$\ $

$\ $

$\ $

\setcounter{section}{8}
{\large\bf 8. The proof of Theorem 1.13} 

$\ $

{\bf 8.1 Notations} In this section $\ncps$ is a fixed 
\setcounter{equation}{0}
non-commutative probability space, such that $\varphi$ is a trace,
and $u, p_{1,1} , p_{1,2} , \ldots , p_{k,1} , p_{k,2} \in \A$
are elements satisfying the conditions (i), (ii), (iii) of 
the Theorem 1.13. It will be handier to prove that the sets
$\{ p_{1,1} \uinv , up_{1,2} \} , \ldots ,
\{ p_{k,1} \uinv , up_{k,2} \}$ are free (this is equivalent 
to the statement of Theorem 1.13, by swapping $p_{j,1}$ with
$p_{j,2}$ for every $1 \leq j \leq k).$ 

We make the notations $p_{j,1} \uinv \egdef a_{j,1},$
$u p_{j,2} \egdef a_{j,2},$ $1 \leq j \leq k.$ In order to
start the proof that the sets $\{ a_{1,1} , a_{1,2} \} ,
\ldots , \{ a_{k,1} , a_{k,2} \}$ are free in $\ncps$,
we will use an abstract nonsense construction: we consider
(and fix for the whole section) another non-commutative
probability space $( \B , \psi ),$ with $\psi$ a trace, 
and elements
$b_{1,1} , b_{1,2} , \ldots , b_{k,1} , b_{k,2} \in \B$ 
such that:

(j) $\mu_{b_{j,1} , b_{j,2}} = \mu_{a_{j,1} , a_{j,2}} ,$
for every $1 \leq j \leq k$ (where the joint distributions
$\mu_{b_{j,1} , b_{j,2}}$ and $\mu_{a_{j,1} , a_{j,2}}$
are considered in $( \B , \psi )$ and $\ncps ,$ respectively).

(jj) the sets
$\{ b_{1,1} , b_{1,2} \} , \ldots , \{ b_{k,1} , b_{k,2} \}$
are free in $( \B , \psi ).$
\newline
Such a construction can of course be done, for instance we can
take $( \B , \psi )$ to be the free product
$ \star_{j=1}^{k} ( \A_{j} , \varphi | \A_{j} ),$ where $\A_{j}$ 
is the unital algebra generated by $\{ a_{j,1} , a_{j,2} \}$ in 
$\A ,$ and then we can take $b_{j,1} , b_{j,2}$ to be just
$a_{j,1} , a_{j,2},$ but viewed in $\B .$ (Note that
$\star_{j=1}^{k} ( \varphi | \A_{j} )$ is a trace, because each 
$\varphi | \A_{j}$ is so, and by Proposition 2.5.3 of \cite{VDN}
- this justifies why we could assume that $\psi$ is a trace.)

In this setting, our goal is to show that
\begin{equation}
\mu_{a_{1,1} , a_{1,2} , \ldots , a_{k,1} , a_{k,2}} \ = \ 
\mu_{b_{1,1} , b_{1,2} , \ldots , b_{k,1} , b_{k,2}} . 
\end{equation}
Indeed, as it is clear from its very definition, the freeness
of a family of subsets can be read from the joint distribution
of the union of those subsets; so in the presence of (8.1),
the freeness of  
$\{ a_{1,1} , a_{1,2} \} , \ldots , \{ a_{k,1} , a_{k,2} \}$ 
is implied by the one of
$\{ b_{1,1} , b_{1,2} \} , \ldots , \{ b_{k,1} , b_{k,2} \} .$ 

$\ $

{\bf 8.2 Remark} For every $1 \leq j \leq k,$ the pair
$( a_{j,1} , a_{j,2} )$ - and hence $( \bjunudoi )$ too - is
$R$-diagonal with determining series
\begin{equation}
f_{j} \ \egdef \ R( \mu_{p_{j,1} \cdot p_{j,2}} ) \
\freestar \  Moeb.
\end{equation}
This follows from Theorem 1.5 and Proposition 1.7 (and the fact
that $a_{j,1} \cdot a_{j,2} = p_{j,1} \cdot p_{j,2} ).$

We have, in other words, that
\begin{equation}
[ R( \mu_{\ajunudoi} )] ( \zjunudoi ) \ = \ 
[ R( \mu_{\bjunudoi} )] ( \zjunudoi ) \ = \ 
f_{j} ( z_{j,1} \cdot z_{j,2} ) +
f_{j} ( z_{j,2} \cdot z_{j,1} ), 
\end{equation}
for every $1 \leq j \leq k.$
The freeness of $\{ \bunuunudoi \} , \ldots , \{ \bkunudoi \}$
enables us to obtain from (8.3) the formula for the $R$-transform
$R( \mu_{\bunuunudoi , \ldots , \bkunudoi} ),$ this is the series
of $2k$ variables
\begin{equation} 
f( \zunuunudoi , \ldots , \zkunudoi ) \ \egdef \ \sum_{j=1}^{k}
f_{j} ( z_{j,1} \cdot z_{j,2} ) + f_{j} ( z_{j,2} \cdot z_{j,1} ).
\end{equation}
(Of course, we don't have the analogue of this fact for the $a$'s
- if we would, the proof would be finished.)

$\ $

A way of re-interpreting (8.2) which will be used later on is
the following:

$\ $

{\bf 8.3 Lemma} For every $1 \leq j \leq k$ and every $n \geq 1,$
the coefficients of $( z_{j,1} \cdot z_{j,2} )^{n},$ 
$( z_{j,2} \cdot z_{j,1} )^{n}$ in  
$[ R( \mu_{\bjunudoi} )] ( \zjunudoi )$ and
$[ R( \mu_{\pjunudoi} )] ( \zjunudoi )$ are (all four of them) equal
to the coefficient of order $n$ in the series $f_{j}$ of (8.2).

$\ $

{\bf Proof} For $R( \mu_{\bjunudoi} )$ this has been explicitly
written in (8.3), while for $R( \mu_{\pjunudoi} )$ we use 
Proposition 5.3, the hypothesis that the pair $( \pjunudoi )$
is diagonally balanced, and the form of $f_{j}$ in (8.2).
{\bf QED}

$\ $

{\bf 8.4 The approach to the proof} of Theorem 1.13 will follow
from now on the same pattern as the one used for the proof of
Theorem 1.5. Indeed, Eqn.(8.1) simply says that for every
$m \geq 1$ and $\lunum \in \{ 1,2 \} ,$ $\hunum \in \kk ,$
we have
\begin{equation}
\varphi ( a_{h_{1},l_{1}} \cdots a_{h_{m},l_{m}} ) \ = \ 
\psi ( b_{h_{1},l_{1}} \cdots b_{h_{m},l_{m}} ) . 
\end{equation}
In the rest of the section we will work on proving (8.5), in a
way which parallels the development of Section 7. We will first
obtain formulas expressing the two sides of (8.5) as summations
over $NC(m).$ From these formulas it will be clear that both 
sides of (8.5) are equal to zero, unless $( \lunum ) \in
\{ 1,2 \} ^{m}$ satisfies a certain balancing condition - in
fact the same as the one stated in 7.5. Starting from that point,
we will assume that the balancing condition holds, and in order 
to complete the proof we will need to throw in a ``geometrical''
argument concerning the circular picture of a non-crossing partition.

$\ $

Both sides of (8.5) are invariant under cyclic permutations of the
monomials involved (because $\varphi$ and $\psi$ are traces).
Thus if we assume in (8.5) that $l_{1} =1,$ we are in fact only 
missing the case when $l_{1} = l_{2} = \cdots = l_{m} =2.$ But
in the latter case we have:

$\ $

{\bf 8.5 Lemma} 
$\varphi ( a_{h_{1},2} a_{h_{2},2} \cdots a_{h_{m},2} )$ =  
$\psi ( b_{h_{1},2} b_{h_{2},2} \cdots b_{h_{m},2} )$ = 0.

$\ $

{\bf Proof} For the $a$'s:  
$\varphi ( a_{h_{1},2} a_{h_{2},2} \cdots a_{h_{m},2} )$ =  
$\varphi ( up_{h_{1},2} up_{h_{2},2} \cdots up_{h_{m},2} )$  
can be written as a coefficient of the series 
$M( \mu_{up_{h_{1},2} , up_{h_{2},2} , \ldots , up_{h_{m},2}} ).$  
But this series is identically zero; indeed, by using Eqn.(3.16)
and the fact that $u$ is free from 
$\{ p_{1,2} , \ldots , p_{k,2} \}$ ), we get:
\[
M( \mu_{up_{h_{1},2} , up_{h_{2},2} , \ldots , up_{h_{m},2}} )
\ = \ M( \mu_{\scriptstyle 
\underbrace{u, \ldots ,u}_{m} } ) \freestar
R( \mu_{p_{h_{1},2} , p_{h_{2},2} , \ldots , p_{h_{m},2}} ).
\]
It is obvious, however, that $M( \mu_{u, \ldots , u} )=0,$ 
and that in general we have $0 \freestar f =0.$

For the $b$'s: $b_{h,2}$ is a part of the $R$-diagonal pair
$( \bhunudoi ),$ and therefore, as remarked in the final paragraph
of Section 6,
we must have $\psi ( b_{h,2}^{n} )=0$ for every $1 \leq h \leq k$
and $n \geq 1.$ But then, if we also take into account that 
$b_{1,2} , b_{2,2} , \ldots , b_{k,2}$ are free, the equality
$\psi ( b_{h_{1},2} b_{h_{2},2} \cdots b_{h_{m},2} )$ = 0
follows directly from the definition of freeness in (1.1).
{\bf QED}

$\ $

{\bf 8.6 Notations} From now on, and until the end of the section,
we fix: $m \geq 1;$ $\lunum \in \{ 1,2 \} ^{m}$ such that $l_{1} =1;$
and $\hunum \in \kk .$ Our goal is to prove (8.5) for this fixed set 
of data.

We denote the $m$-tuple $( \lunum ) \in \{ 1,2 \} ^{m}$ by $\ee ,$
and we denote by $n$ $( 1 \leq n \leq m)$ the number of occurrences 
of 1 in $\ee .$ 

We will use the geometrical objects constructed in Sections 7.1, 7.2
above. That is, we consider again the circle of radius 1 and the 
points $\Punum , \Qunum$ sitting on it, and positioned in the way
described in 7.1; and we consider again the complementation map 
$\cq : NC(m) \rightarrow NC(m)$ described in 7.2.

$\ $

{\bf 8.7 Remark} Consider the product 
$a_{h_{1},l_{1}} a_{h_{2},l_{2}} \cdots a_{h_{m},l_{m}}$
appearing in the right-hand side of (8.5). If in this product
we write back each $a_{h_{i},1}$ as $p_{h_{i},1} \uinv$ and
each $a_{h_{i},2}$ as $u p_{h_{i},2}$, we obtain an expression
(monomial) of length $2m$ in $\punuunudoi , \ldots ,
\pkunudoi , u, \uinv .$ By looking just at the $p$'s in the
latter monomial, we see that they are
$p_{h_{1},l_{1}} , p_{h_{2},l_{2}} , \ldots , p_{h_{m},l_{m}} ,$
exactly in this order, and placed on a certain set of positions
$I \subseteq \mmmm ,$ with $|I|=m.$ On the complementary set of 
positions $J = \mmmm \setminus I$ of our monomial of length $2m$
we have factors of $u$ and $\uinv ;$ we denote them as
$u^{\lambda_{1}}, u^{\lambda_{2}}, \ldots ,u^{\lambda_{m}},$ in
the order in which they appear from left to right. It is clear that
$\lambda_{i}, \ 1 \leq i \leq m,$ is entirely determined by 
$l_{i},$ via the formula:
$l_{i} =1 \Rightarrow \lambda_{i} = -1,$
$l_{i} =2 \Rightarrow \lambda_{i} = +1.$

We can thus say that the product
$a_{h_{1},l_{1}} a_{h_{2},l_{2}} \cdots a_{h_{m},l_{m}}$
is obtained by shuffling together 
$p_{h_{1},l_{1}} p_{h_{2},l_{2}} \cdots p_{h_{m},l_{m}}$
and $u^{\lambda_{1}} u^{\lambda_{2}} \cdots u^{\lambda_{m}},$ 
where the $p$'s have to sit on the positions indicated by $I$,
and the $u$'s have to sit on the positions indicated by 
$J = \mmmm \setminus I.$ Note that $I$ and $J$ are exactly the 
same as the ones appearing in the proof of Proposition 7.3.

$\ $

{\bf 8.8 Proposition} In the notations established in 8.1,
8.6, 8.7, we have:
\[
\varphi ( a_{h_{1},l_{1}} a_{h_{2},l_{2}} \cdots 
a_{h_{m},l_{m}} ) \ =
\]
\begin{equation}
= \ \sum_{\sigma \in NC(m)}
[ \mbox{coef } ( \hlunum ); \sigma ]
( R( \mu_{\punuunudoi , \ldots , \pkunudoi} )) \cdot
\end{equation}
\[
[ \mbox{coef } ( \llunum ); \cq ( \sigma ) ]
( M( \mu_{u, \uinv} ))
\]
and
\[
\psi ( b_{h_{1},l_{1}} b_{h_{2},l_{2}} \cdots 
b_{h_{m},l_{m}} ) \ =
\]
\begin{equation}
= \ \sum_{\sigma \in NC(m)}
[ \mbox{coef } ( \hlunum ); \sigma ]
( R( \mu_{\bunuunudoi , \ldots , \bkunudoi} )) .
\end{equation}
(In (8.6), (8.7) the series 
$R( \mu_{\punuunudoi , \ldots , \pkunudoi} )$ and
$R( \mu_{\bunuunudoi , \ldots , \bkunudoi} )$
are acting in the variables
$\zunuunudoi , \ldots , \zkunudoi$ - same as in (8.4),
for instance. The series $M( \mu_{u, \uinv} )$ is viewed as
acting in the two variables $z_{+1}$ and $z_{-1}.)$

$\ $

{\bf Proof} Equation (8.7) is just the moment-cumulant formula,
see Section 3.5 above.

The argument proving (8.6) is very similar to the one used in 
the proof of Proposition 7.3, and for this reason we will only 
outline its main steps. We view 
$\varphi ( a_{h_{1},l_{1}} a_{h_{2},l_{2}} \cdots 
a_{h_{m},l_{m}} )$ as a coefficient of length $2m$ of the
series $M( \mu_{\punuunudoi , \ldots , \pkunudoi, u, \uinv} ),$
and we then expand it as a summation over $NC(2m)$ by using 
the moment-cumulant formula (i.e., the appropriate version of 
Eqn.(3.12) in 3.5). Due to the fact that
$\{ \punuunudoi , \ldots , \pkunudoi \}$ is free from 
$\{ u, \uinv \} ,$ what we arrive to is the formula (paralleling
(7.4) in the proof of 7.3):
\[
\varphi ( a_{h_{1},l_{1}} a_{h_{2},l_{2}} \cdots 
a_{h_{m},l_{m}} ) \ =
\]
\begin{equation}
= \ \sum_{ \begin{array}{c}
{\scriptstyle \sigma , \tau \in NC(m)}  \\
{\scriptstyle I,J - compatible} 
\end{array} }
[ \mbox{coef } ( \hlunum ); \sigma ]
( R( \mu_{\punuunudoi , \ldots , \pkunudoi} )) \cdot
\end{equation}
\[
[ \mbox{coef } ( \llunum ); \tau ) ] ( R( \mu_{u, \uinv} )).
\]
Then we process the right-hand side of (8.8) exactly in the
way the right-hand side of (7.4) was processed in the proof 
of Proposition 7.3, and this leads to the right-hand side of 
(8.6). { \bf QED}

$\ $

{\bf 8.9 Corollary} Let $\ee = ( \lunum )$ be the string of
1's and 2's of Notations 8.6. If the number $n$ of occurrences
of 1 in $\ee$ does not equal the number $m-n$ of occurrences of
2 in $\ee ,$ then 
$\varphi ( a_{h_{1},l_{1}} a_{h_{2},l_{2}} \cdots 
a_{h_{m},l_{m}} )$ =
$\psi ( b_{h_{1},l_{1}} b_{h_{2},l_{2}} \cdots 
b_{h_{m},l_{m}} )$ = 0.

$\ $

{\bf Proof} For the $a$'s: Due to how $\llunum$ are determined 
by $\lunum$ (see Remark 8.7), we have in this case that the 
number $n$ of $(-1)$'s in $( \llunum )$ is different from the 
number $m-n$ of $1$'s in $( \llunum )$. But then it is
immediate that 
$[ \mbox{coef } ( \llunum ); \tau] ( M( \mu_{u, \uinv} ))$
= 0 for every $\tau \in NC(m),$ and hence all the terms in
the sum of (8.6) are vanishing.

For the $b$'s: in (8.4) we have an explicit formula for
$R( \mu_{\bunuunudoi , \ldots , \bkunudoi} ),$ the inspection
of which makes clear (in view of the hypothesis of the present 
corollary) that all the terms of the sum in (8.7) are vanishing.
{\bf QED}

$\ $

We therefore make, exactly as we did in 7.5, the following

$\ $

{\bf 8.10 Supplementary hypothesis:} From now on we will assume
that the string $\ee = ( \lunum )$ of 8.6 contains an equal
number of 1's and 2's (i.e., $m$ is even and $n=m/2).$

$\ $

Moreover, we will consider the notion of an {\em $\ee$-alternating}
partition in $NC(m),$ which is exactly the one defined in 7.6. 
Recall that the set of all the $\ee$-alternating partitions in
$NC(m)$ is denoted by $\eealt .$

$\ $

{\bf 8.11 Proposition} A partition $\sigma \in NC(m)$ is
$\ee$-alternating if and only if it has the following properties:

$1^{o}$ every block of $\sigma$ has at least two elements;

$2^{o}$ there exists no block $B$ of $\sigma$ such that: $|B|$
is odd, and $( \lunum )|B$ is a cyclic permutation of
$( \underbrace{1,2, \ldots ,1,2,1}_{|B|} )$ or of
$( \underbrace{2,1, \ldots ,2,1,2}_{|B|} )$
(where the notation $( \lunum )|B$ is in the sense of 3.2.$2^{o});$

$3^{o}$ every block of the complementary partition $\cq ( \sigma )$
has an even number of elements.

$\ $

{\bf Proof} ``$\Rightarrow$'' As remarked immediately after the
Definition 7.6, every block of $\sigma$ has an even number of
elements (this implies $1^{o}$ and $2^{o}).$ The same holds for
$\cq ( \sigma ),$ because $\cq ( \sigma )$ is also in $\eealt ,$
by Proposition 7.7.

``$\Leftarrow$'' By contradiction, we assume that $\sigma$ is not
$\ee$-alternating. This means that we can find 
$1 \leq i<j \leq m$ such that $i$ and $j$ are in the same 
block of $\sigma ,$ and consecutive in that block (in the sense of
2.2), and such that moreover $l_{i} =l_{j}.$ From all
the pairs $(i,j)$ with these properties, 
we choose one, $( \iz , \jz ),$ such that the length 
of the segment $P_{\iz}P_{\jz}$ is minimal. (We are using in 
this proof the circular picture involving the points 
$\Punum , \Qunum$ introduced in 7.1, 7.2.)

We denote by $\sps$ the open half-plane that is determined by the
line through $P_{\iz}$ and $P_{\jz}$, and that does not contain the 
center of the circle. (If the line $P_{\iz}P_{\jz}$ goes through the 
center of the circle, any of the two open half-planes determined 
by it can be chosen as $\sps .)$ The minimality assumption on 
$( \iz , \jz )$ implies that:
\begin{equation}
\left\{  \begin{array}{l}
{ \mbox{ if $1 \leq i < j \leq m$ are in the same block of
$\sigma ,$ and consecutive in that block, }  }  \\
{ \mbox{ and if in addition we have that $P_{i},P_{j} \in
\sps ,$ then necessarily $l_{i} \neq l_{j}$}  }  
\end{array}   \right.
\end{equation}
(this is simply because the length of $P_{i}P_{j}$ is strictly
smaller than the one of $P_{\iz}P_{\jz}).$ Note that this argument
gives in fact that $l_{i} \neq l_{j}$ even if we only assume
that $P_{i},$ $P_{j}$ lie in the closure of $\sps$, but
$(i,j) \neq ( \iz , \jz ).$

Let $B$ be the block of $\sigma$ containing $\iz$ and $\jz .$ 
We claim that $\{ P_{j} \ | \ j \in B \} \cap \sps = \emptyset .$ 
Indeed, this is obvious if $B$ is reduced to $\{ \iz , \jz \}$, 
so let us assume that $|B| \geq 3.$ Since $\iz , \jz$ 
are consecutive in $B,$ we have that $\{ P_{j} \ | \ j \in B \}$ 
is contained in one of the closed half-planes determined by 
$P_{\iz}P_{\jz} ;$ so if we would assume 
$\{ P_{j} \ | \ j \in B \} \cap \sps \neq \emptyset ,$ 
then it would follow that 
$\{ P_{j} \ | \ j \in B , \ j \neq \iz , \jz \} \subseteq \sps .$ 
But then the remark concluding the preceding paragraph would
imply that $l_{i} \neq l_{j}$ whenever $i<j$ are consecutive in
$B$, with the exception of the case when $(i,j) = ( \iz , \jz ),$
and this would violate property $2^{o}$ of $\sigma .$ 

Now, the set $\{ j \ | \ 1 \leq j \leq m, \ P_{j} \in \sps \}$
is either void or a union of blocks of $\sigma ,$ because of the
non-crossing character of $\sigma$ (and because of what was proved 
about the block $B \ni \iz , \jz ).$ Remark
that for every block $B'$ of $\sigma$ which is contained in
$\{ j \ | \ 1 \leq j \leq m, \ P_{j} \in \sps \} ,$ we must have
that $|B'|$ is even; indeed, from (8.9) it follows that 
$l_{i} \neq l_{j}$ for every $i,j \in B'$ which are consecutive 
in $B' ,$ and this couldn't happen if $|B'|$ would be odd. 
Consequently, the set
$\{ j \ | \ 1 \leq j \leq m, \ P_{j} \in \sps \}$
has an even cardinality (possibly zero).

By recalling how the points $\Qunum$ were constructed, we next 
note that $\{ j \ | \ 1 \leq j \leq m, \ P_{j} \in \sps \}$
$\subset$
$\{ j \ | \ 1 \leq j \leq m, \ Q_{j} \in \sps \} ,$ and moreover
that the set-difference $\{ j \ | \ 1 \leq j \leq m,$
$Q_{j} \in \sps , \ P_{j} \not\in \sps \}$ has exactly one element,
which is either $\iz$ or $\jz .$ The latter assertion amounts to the 
following two facts, both obvious: (a) if $1 \leq j \leq m$ is such
that $P_{j}$ does not belong to the closure of $\sps ,$ then 
$Q_{j} \not\in \sps ;$ and (b) $Q_{\iz}$ and $Q_{\jz}$ lie in 
opposite open half-planes determined by the line $P_{\iz} P_{\jz},$
hence exactly one of them is in $\sps$ (the assumption
$l_{\iz} = l_{\jz}$ is of course crucial for (b)). The conclusion
of this paragraph is that the cardinality
$| \{ j \ | \ 1 \leq j \leq m, \ Q_{j} \in \sps \} |$ = 
1 + $| \{ j \ | \ 1 \leq j \leq m, \ P_{j} \in \sps \} |$
is an odd number.

But because of how the definition of the complementation map 
$\cq$ was made in 7.2, the set 
$\{ j \ | \ 1 \leq j \leq m, \ Q_{j} \in \sps \}$ is a union of
blocks of $\cq ( \sigma ).$ Hence, in view of the property $3^{o}$
of $\sigma ,$ 
$| \{ j \ | \ 1 \leq j \leq m, \ Q_{j} \in \sps \} |$ is an even
number - contradiction. {\bf QED}

$\ $

{\bf 8.12 The conclusion of the proof} Recall that what we need to 
prove is the equality (8.5) of 8.4, for the $m \geq 1,$
$\lunum \in \{ 1,2 \} ,$ $\hunum \in \kk$ that were fixed in 8.6,
and where $\ee = ( \lunum )$ satisfies the additional hypothesis
8.10.

Let us express both quantities 
$\varphi ( a_{h_{1},l_{1}} \cdots a_{h_{m},l_{m}} )$  and
$\psi ( b_{h_{1},l_{1}} \cdots b_{h_{m},l_{m}} )$ 
appearing in (8.5) as summations over $NC(m),$
in the way shown in Proposition 8.8 (Eqns.(8.6) and (8.7),
respectively). We will prove that for every $\sigma \in NC(m),$
the terms indexed by $\sigma$ in the two sums of (8.6) and (8.7)
coincide - this will of course imply that the sums are equal.

Up to now, the $m$-tuple $( \hunum ) \in  \kk^{m} ,$
which is part of our data, didn't play any role. Let us view this
$m$-tuple as a function from $\mm$ to $\kk ,$ and denote by 
$\Lunuk$ its level sets; i.e., $L_{h} \egdef \{ i \ |$
$1 \leq i \leq m, \ h_{i} =h \} ,$ for $1 \leq h \leq k.$
We will say about a partition $\sigma \in NC(m)$ that it is
{\em $h$-acceptable} if every block of $\sigma$ is contained
in one of the level sets $\Lunuk ;$ and we will denote by
$\hacc (m)$ the set of all $h$-acceptable partitions in $NC(m).$

The general term of the sum in (8.7) is a product of coefficients
of $R( \mu_{\bunuunudoi , \ldots , \bkunudoi} ),$ made in the way
dictated by $\sigma \in NC(m)$ (according to the recipe of
3.2.$3^{o}).$ By using the freeness of the sets
$\{ \bunuunudoi \} , \ldots , \{ \bkunudoi \}$ and Theorem 3.6,
it is easily seen that this product is zero whenever $\sigma$ is 
not $h$-acceptable; while for $\sigma \in \hacc (m)$ we get the
formula:
\begin{equation}
[ \mbox{coef } ( \hlunum ); \sigma ]
( R( \mu_{\bunuunudoi , \ldots , \bkunudoi} )) \ =
\end{equation}
\[
= \ \prod_{ \begin{array}{c}
{\scriptstyle 1 \leq h \leq k}  \\   
{\scriptstyle such \ that}  \\   
{\scriptstyle L_{h} \neq \emptyset}     
\end{array}  }  \ \{
\ \prod_{ \begin{array}{c}
{\scriptstyle B \ block \ of \ \sigma}  \\   
{\scriptstyle such \ that}  \\   
{\scriptstyle B \subseteq L_{h}}     
\end{array}  }   \
[ \mbox{coef } ( \hlunum ) | B]
( R( \mu_{\bhunudoi} )) \ \} ,
\]
(where the notation for the restricted $|B|$-tuple $( \hlunum ) | B$
is taken from 3.2.$2^{o}).$ We next recall from 8.2 that the pair 
$( \bhunudoi )$ is $R$-diagonal, for every $1 \leq h \leq k;$ hence
even if we assume that $\sigma \in \hacc (m),$ the particular form
of $R( \mu_{\bhunudoi} )$ will force the product in (8.10) to still
be equal to zero, unless we also assume that the $|B|$-tuple
$( \lunum )|B$ is $(1,2,1,2, \ldots ,1,2)$ or $(2,1,2,1, \ldots ,
2,1)$ for every block $B$ of $\sigma .$ In other words, if we want
the contribution of $\sigma \in NC(m)$ to the sum (8.7) to be 
non-zero, then we must also impose (besides the condition
$\sigma \in \hacc (m))$ that $\sigma$ is $\ee$-alternating.

To conclude the discussion concerning (8.7), let us consider
$\sigma \in \hacc (m) \cap \eealt .$ Then the product
in (8.10) can be prelucrated by using Lemma 8.3, and this yields
the formula
\[
[ \mbox{coef } ( \hlunum ); \sigma ]
( R( \mu_{\bunuunudoi , \ldots , \bkunudoi} )) \ =
\]
\begin{equation}
= \ \prod_{ \begin{array}{c}
{\scriptstyle 1 \leq h \leq k}  \\   
{\scriptstyle such \ that}  \\   
{\scriptstyle L_{h} \neq \emptyset}     
\end{array}  }  \ \{
\ \prod_{ \begin{array}{c}
{\scriptstyle B \ block \ of \ \sigma}  \\   
{\scriptstyle such \ that}  \\   
{\scriptstyle B \subseteq L_{h}}     
\end{array}  }   \
[ \mbox{coef } ( |B|/2 ) ] ( f_{h} ) \ \} 
\end{equation}
(where $\funuk$ are the series of one variable discussed in
Remark 8.2).

\vspace{10pt}

We now start looking at the term indexed by $\sigma$ in the sum
(8.6). From the parallel discussion made above, we see that 
we have to consider three cases: (a) $\sigma \in NC(m) \setminus
\hacc (m);$ (b) $\sigma \in \hacc (m) \setminus \eealt ;$
and (c) $\sigma \in \hacc (m) \cap \eealt .$ More precisely,
we have to show that the term indexed by $\sigma$ in (8.6) is
zero in the cases (a) and (b), and is given by the right-hand
side of (8.11) in the case (c).

\vspace{6pt}

Case (a) is easy. Indeed, by using the freeness of 
$\{ \punuunudoi \} , \ldots , \{ \pkunudoi \}$ and Theorem 3.6,
it is easily seen that in this case 
$[ \mbox{coef } ( \hlunum ); \sigma ]$
$( R( \mu_{\punuunudoi , \ldots , \pkunudoi} ))$ equals zero.

\vspace{6pt}

Case (c) is also easy. Indeed, in this case the $h$-acceptability 
of $\sigma$ implies that 
\newline
$[ \mbox{coef } ( (h_{1},l_{1}),$
$\ldots , (h_{m},l_{m}) ); \sigma ]$
$( R( \mu_{\punuunudoi , \ldots , \pkunudoi} ))$ has an expression of
the kind shown in the right-hand side of (8.10) (with the $b$'s
replaced by $p$'s); and after that, due to the assumption that
$\sigma \in \eealt ,$ this expression can be prelucrated by using 
Lemma 8.3, and yields exactly the right-hand side of (8.11).
On the other hand, in case (c) we also have $\cq ( \sigma ) \in
\eealt $ (by Proposition 7.7); by using this fact, and by following
how $( \llunum ) $ is determined by $\ee = ( \lunum )$ (see Remark
8.7), one sees immediately that we have
$[ \mbox{coef } ( \lambda_{1}, \ldots , \lambda_{m} ); \cq ( \sigma )]$ 
$(M ( \mu_{u, \uinv} ))$ =1.

\vspace{6pt}

Finally, let us consider the case (b). The $h$-acceptability of 
$\sigma$ implies that (similarly to (8.10)) we have the formula
\[
[ \mbox{coef } ( \hlunum ); \sigma ]
( R( \mu_{\punuunudoi , \ldots , \pkunudoi} )) \ =
\]
\begin{equation}
= \ \prod_{ \begin{array}{c}
{\scriptstyle 1 \leq h \leq k}  \\   
{\scriptstyle such \ that}  \\   
{\scriptstyle L_{h} \neq \emptyset}     
\end{array}  }  \ \{
\ \prod_{ \begin{array}{c}
{\scriptstyle B \ block \ of \ \sigma}  \\   
{\scriptstyle such \ that}  \\   
{\scriptstyle B \subseteq L_{h}}     
\end{array}  }   \
[ \mbox{coef } ( \hlunum ) | B]
( R( \mu_{\phunudoi} )) \ \} .
\end{equation}
We have to show that either the product in 
(8.12) is zero, or $[ \mbox{coef } ( \llunum ); \cq ( \sigma ) ]$
\newline
$(M ( \mu_{u, \uinv} ))$ is zero. We will obtain this from the hypothesis 
that $\sigma \not\in \eealt ,$ in the equivalent formulation that
comes out by negating Proposition 8.11. That is, we know that because
$\sigma \not\in \eealt ,$ one of the following three things must
happen:

$1^{o}$ either: $\sigma$ has a block $B$ with one element. In this
case the product in (8.12) contains a factor which is a coefficient
of length one of the series $R( \mu_{\phunudoi} ),$
$1 \leq h \leq k;$ but any such coefficient is zero, because the 
pairs $( \phunudoi )$ are diagonally balanced.

$2^{o}$ or: $\sigma$ has a block $B$ such that $|B|$ is odd and
$( \lunum )|B$ is a cyclic permutation of $(1,2, \ldots ,1,2,1)$ 
or of $(2,1, \ldots ,2,1,2)$. Let $h \in \kk$
be such that $B$ is contained in the level set $L_{h}.$ Then 
$[ \mbox{coef } ( (h_{1},l_{1}), \ldots , (h_{m},l_{m}) )|B ] 
( R( \mu_{\phunudoi} ))$ (which is the same thing as
$[ \mbox{coef } ( (h,l_{1} ), \ldots ,(h, l_{m} ))|B ] 
( R( \mu_{\phunudoi} ))$, by the definition of $L_{h})$, is zero
- because $( \phunudoi )$ is diagonally balanced, and by Remark 5.2.

$3^{o}$ or: $\cq (\sigma )$ has a block containing an odd number of 
elements. But then it is clear that 
$[ \mbox{coef } ( \llunum ); \cq ( \sigma ) ]
(M ( \mu_{u, \uinv} ))$ =0. This concludes the discussion of case
(b), and the proof. {\bf QED}

$\ $

$\ $

$\ $


\begin{thebibliography}{99}

\bibitem{Ba} T. Banica.
On the polar decomposition of circular variables, preprint.

\bibitem{Bi} Ph. Biane.
Some properties of crossings and partitions, preprint.

\bibitem{DRS} P. Doubilet, G.-C. Rota, R. Stanley.
On the foundations of the combinatorial theory (VI):
the idea of generating function, in Proc. of the 6th
Berkeley symposium on mathematical statistics and probability,
Lucien M. Le Cam et al. editors, University of California Press,
1972, 267-318.

\bibitem{D1} K. Dykema.
Interpolated free group factors, Pacific J. Math. 163(1994),
123-135.

\bibitem{D2} K. Dykema.
Free products of hyperfinite von Neumann algebras and free
dimension, Duke J. Math. 69(1993), 97-119. 

\bibitem{E} P. Edelman.
Chain enumeration and non-crossing partitions, Discrete
Math. 31(1980), 171-180.

\bibitem{K} G. Kreweras.
Sur les partitions non-croisees d'un cycle, Discrete Math.
1(1972), 333-350.

\bibitem{N} A. Nica.
$R$-transforms of free joint distributions, and non-crossing 
partitions, to appear in the Journal of Functional Analysis.

\bibitem{NS1} A. Nica, R. Speicher.
A ``Fourier transform'' for multiplicative functions on non-crossing
partitions, preprint.

\bibitem{NS2} A. Nica, R. Speicher.
On the multiplication of free $n$-tuples of non-commutative random
variables, preprint.

\bibitem{R1} F. Radulescu.
The fundamental group of the von Neumann algebra of a free group
with infinitely many generators is $\R_{+} \setminus \{ 0 \} ,$
J. Amer. Math. Soc. 5(1992), 517-532.

\bibitem{R2} F. Radulescu.
Random matrices, amalgamated free products and subfactors of the
von Neumann algebra of a free group, of noninteger index,
Inventiones Math. 115(1994), 347-389.

\bibitem{SU} R. Simion, D. Ullman.
On the structure of the lattice of non-crossing partitions,
Discrete Math. 98(1991), 193-206.

\bibitem{S1} R. Speicher.
Multiplicative functions on the lattice of non-crossing partitions
and free convolution, Math. Annalen 298(1994), 611-628.

\bibitem{S2} R. Speicher.
Combinatorial theory of the free product with amalgamation and
operator-valued free probability theory, Habilitationsschrift,
Heidelberg 1994.

\bibitem{V1} D. Voiculescu.
Symmetries of some reduced free product C*-algebras, in
Operator algebras and their connection with topology and ergodic
theory, H. Araki et al. editors
(Springer Lecture Notes in Mathematics, volume 1132, 1985), 556-588.

\bibitem{V2} D. Voiculescu.
Addition of certain non-commuting random variables, J. Functional
Analysis 66(1986), 323-346.

\bibitem{V3} D. Voiculescu.
Multiplication of certain non-commuting random variables,
J. Operator Theory 18(1987), 223-235.

\bibitem{V4} D. Voiculescu.
Limit laws for random matrices and free products, Inventiones
Math. 104(1991), 201-220.

\bibitem{V5} D. Voiculescu.
Circular and semicircular systems and free product factors, in
Operator algebras, unitary representations, enveloping algebras,
and invariant theory, A. Connes et al. editors, Birkh\"{a}user, 
1990, 45-60.

\bibitem{VDN} D. Voiculescu, K. Dykema, A. Nica.
Free random variables, CRM Monograph Series, volume 1, AMS, 1992.

\end{thebibliography}
\end{document}